\documentclass[10pt]{article}\textwidth 160mm\textheight 235mm
\usepackage{amsfonts} \usepackage{amssymb} \usepackage{graphicx}
\usepackage{amsmath} \usepackage{amsthm} \usepackage{mathrsfs}
\usepackage{tikz} \usepackage{cancel} \usepackage{todonotes}
\usetikzlibrary{matrix} \usetikzlibrary{arrows,shapes}
\usepackage{cite} \usepackage{hyperref}
\usepackage{amsmath}
\usepackage{amsthm}
\usepackage{amsfonts}
\usepackage{multirow}
\usepackage{amssymb}
\usepackage{graphicx}
\usepackage{tcolorbox}
\usepackage{booktabs}
\usepackage{tikz-cd}

\newtheorem{Theorem}{Theorem}[section]
\newtheorem{Proposition}[Theorem]{Proposition}
\newtheorem{Remark}[Theorem]{Remark}
\newtheorem{Lemma}[Theorem]{Lemma}
\newtheorem{Corollary}[Theorem]{Corollary}
\newtheorem{Definition}[Theorem]{Definition}
\newtheorem{Example}[Theorem]{Example}

\usepackage{hyperref}
\hypersetup{colorlinks=true,urlcolor=blue}
\expandafter\let\expandafter\oldproof\csname\string\proof\endcsname
\let\oldendproof\endproof
\renewenvironment{proof}[1][\proofname]{
\oldproof[\ttfamily\scshape \bf #1.]
}{\oldendproof}

\DeclareMathOperator*{\argmin}{arg\,min}
\def\tilde{\widetilde}
\def\emp{\emptyset}  
\def\dom{{\rm dom}\,}  \def\epi{{\rm epi\,}}

\def\ox{\bar{x}} \def\oy{\bar{y}} 
  \def\disp{\displaystyle}
\def\tto{\rightrightarrows} \def\Hat{\widehat}

\def\ra{\rangle}
\def\la{\langle}
\def\ve{\varepsilon}
\def\e{\varepsilon}
 \def\ov{\bar{v}}

\def\gph{\mbox{\rm gph}\,} \def\epi{\mbox{\rm epi}\,}
 \def\dom{\mbox{\rm dom}\,}

 \def\dn{\downarrow} \def\O{\Omega}
\def\ph{\varphi} \def\emp{\emptyset} \def\st{\stackrel}
\def\oR{\overline{\R}} 
\def\lm{\lambda} 
\def\gg{\gamma} \def\dd{\delta}
\def\al{\alpha}  \def\N{{\rm I\!N}}
\def\R{{\rm I\!R}}
 \def\vt{\vartheta} \def\e{\varepsilon} 
\def\Limsup{\mathop{{\rm Lim}\,{\rm sup}}}
\def\Liminf{\mathop{{\rm Lim}\,{\rm inf}}}

\def\Limsup{\mathop{{\rm Lim}\,{\rm sup}}}

\numberwithin{equation}{section}

\title{\bf
Variational and Strong Variational Convexity\\ in Infinite-Dimensional Variational Analysis}
\author{Pham Duy Khanh\footnote{Department of Mathematics, Ho Chi Minh City University of Education, Ho Chi Minh City, Vietnam. E-mail: pdkhanh182@gmail.com} \quad Vu Vinh Huy Khoa\footnote{Department of Mathematics, Wayne State University, Detroit, Michigan, USA. E-mail: khoavu@wayne.edu. Research of this author was partly supported by the US National Science Foundation under grant DMS-2204519.}\quad Boris S. Mordukhovich\footnote{Department of Mathematics, Wayne State University, Detroit, Michigan, USA. E-mail: aa1086@wayne.edu. Research of this author was partly supported by the US National Science Foundation under grants DMS-1808978 and DMS-2204519, by the Australian Research Council under Discovery Project DP-190100555, and by Project~111 of China under grant D21024.}\quad Vo Thanh Phat\footnote{Department of Mathematics, Wayne State University, Detroit, Michigan, USA. E-mail: phatvt@wayne.edu. Research of this author was partly supported by the US National Science Foundation under grants DMS-1808978 and DMS-2204519.}}
\begin{document}
\maketitle\vspace*{-0.2in}
\noindent
{\small{\bf Abstract}. This paper is devoted to a systematic study and characterizations of the fundamental notions of variational and strong variational convexity for lower semicontinuous functions. While these notions have been quite recently introduced by Rockafellar, the importance of them has been already recognized and documented in finite-dimensional variational analysis and optimization. Here we address general infinite-dimensional settings and derive comprehensive characterizations of both variational and strong variational convexity notions by developing novel techniques, which are essentially different from finite-dimensional counterparts. As a consequence of the obtained characterizations, we establish new quantitative and qualitative relationships between strong variational convexity and tilt stability of local minimizers in appropriate frameworks of Banach spaces.\\[1ex]
{\bf Key words}. variational analysis and optimization, variational and strong variational convexity, local maximal monotonicity, proximal mappings and Moreau envelopes, tilt-stable local minimizers, generalized differentiation\\[1ex]
{\bf Mathematics Subject Classification (2020)} 49J52, 49J53, 47H05, 90C30, 90C45\vspace*{-0.15in}}

\section{Introduction}\label{intro}\vspace*{-0.05in}

The notions of {\em variational} and {\em strong variational} convexity of lower semicontinuous functions have been recently introduced by Rockafellar \cite{rtr251} in finite-dimensional spaces, and then they have been applied by him to the design and justification of numerical algorithms of the proximal and augmented Lagrangian types for solving some classes of nonsmooth optimization problems; see \cite{roc,r22}. The motivation behind these notions was the fact that just certain {\em local} (vs.\ global) {\em monotonicity} properties of (limiting) {\em subgradient mappings} were needed for developing numerical algorithms; cf.\ also Pennanen \cite{Pen02}. While it has been well recognized in variational analysis that the crucial {\em global} maximal monotonicity of subgradient mappings is {\em equivalent} to convexity, {\em local} maximal monotonicity opens the door to introducing the new class of {\em variationally convex} functions, and similarly for the strong counterparts. 

On the other lines of research, the usefulness of variational convexity has been recently revealed in the developments of the {\em generalized Newton methods} \cite{BorisKhanhPhat,kmptmp,mor-sar} to solve subgradient inclusions and optimization problems generated by a large class of {\em prox-regular} functions via reducing them to problems of ${\cal C}^{1,1}$ optimization by using {\em Moreau envelopes}. It is shown in \cite{kmp22convex} that the desired convexity (strong convexity) of Moreau envelopes is {\em equivalent} to the variational convexity (strong variational convexity) of the corresponding cost functions. We also mention remarkable applications of variational convexity to {\em variational sufficiency} in constrained optimization \cite{kmp22convex,rtr251,roc,ding} with its algorithmic implementations.\vspace*{0.03in}

All the aforementioned developments concern problems in finite-dimensional spaces with the essential usage of finite-dimensional geometry in formulations and proofs of the obtained results. For its own sake and having in mind applications to infinite-dimensional optimization problems, including dynamic optimization and optimal control of ODE and PDE systems, we now aim at developing variational convexity theory in appropriate settings of {\em Banach spaces}. This is done below by using powerful tools of infinite-dimensional variational analysis and generalized differentiation, epigraphical convergence of extended-real-valued functions, and Banach space geometry. The major results obtained in this way include {\em complete characterizations} of (strong) variational convexity via local (strong) maximal monotonicity of subgradient mappings, via (strong) convexity of tilted and usual Moreau envelopes, relationships between strong variational convexity and tilt stability of local minimizers, etc. The reader will see below that different characterizations require appropriate geometric assumptions on the Banach spaces in question, which all hold in the case of {\em Hilbert spaces}.\vspace*{0.03in}

The rest of the material is organized as follows.          Section~\ref{sec:global-Geometry} contains those definitions and {\em preliminaries}, which are broadly used throughout the entire paper. In Section~\ref{str-conv}, we discuss several notions of {\em strong convexity} for extended-real-valued functions and {\em strong monotonicity} for set-valued mappings with highlighting relationships between them in the case of subgradient operators. Section~\ref{sec:localmonoproperties} concerns {\em local monotonicity} of set-valued mappings together with its strong counterpart and presents the {\em resolvent characterization} (of a local Minty type) of local maximal monotonicity in reflexive Banach spaces.

In Section~\ref{sec:localsub}, we define both notions of variational and strong variational convexity for functions on Banach spaces and derive explicit descriptions of their {\em subdifferential graphs}. This result is instrumental to establish in Section~\ref{sec:vc-Moreau-localmax} our major characterizations of {\em variational convexity} via {\em local maximal monotonicity of subgradient mappings} and {\em local convexity of Moreau envelopes} under various geometric assumptions on the Banach spaces in question. In this section, we also derive {\em quantitative} characterizations of {\em strong variational convexity} with establishing relationships between the corresponding moduli.

Section~\ref{sec:2nd} involves tools of {\em generalized differentiation} in variational analysis to study variational and strong variational convexity. Our characterizations here are obtained in Hilbert spaces via {\em second-order subdifferentials}.  In Section~\ref{sec:tilt}, we establish relationships between {\em strong variational convexity} and {\em tilt stability} of local minimizers for a general class of optimization problems. The {\em concluding} Section~\ref{sec:conc} summarizes the obtained results and discusses some directions of our future research. 
\vspace*{-0.15in}

\section{Preliminaries from Variational Analysis}\label{sec:global-Geometry}\vspace*{-0.05in}

In this section, we define and discuss some notions from variational analysis and generalized differentiation that are frequently used in what follows. More details and references can be found in the books \cite{Mordukhovich18,Rockafellar98} and \cite{Mordukhovich06,Thibault} in finite and infinite dimensions, respectively. 

Unless otherwise stated, all the spaces under consideration are assumed to be {\em Banach} with the generic symbols $\|\cdot\|$ for norms and $\langle\cdot,\cdot\rangle$ for the canonical pairing between a space $X$ and its topological dual $X^*$. We use `$\rightarrow$', `$\overset{w}{\rightarrow}$', and `$\overset{w^*}{\rightarrow}$' to indicate the strong/norm, weak, and weak$^*$ convergence, respectively. Product spaces $X\times Y$ are normed by $\|(x,y)\|:= \|x\|+\|y\|$ whenever $x\in X$ and $y\in Y$. The symbols $\mathbb{B}_{r}(x)$ and $B_r (x)$, signify, respectively, the closed and open balls centered at $x$ with radius $r>0$. Throughout the paper, we consider set-valued mappings/multifunctions $F\colon X\tto Y$ with the {\em domain} $\dom F:=\{x\in X\;|\;F(x)\ne\emp\}$ and the {\em graph} $\gph F:=\{(x,y)\in X\times Y\;|\;y\in F(x)\}$. As usual, $\N:=\{1,2,\ldots\}$.

The (Painlev\'e-Kuratowski) {\em sequential outer limit} of a set-valued  mapping $F:X\rightrightarrows X^*$ as $x\rightarrow \ox$ is
\begin{eqnarray}\label{eq:P-K}
\Limsup\limits_{x \rightarrow \ox} F(x):=\big\{ x^*\in X^* \;\big|\; \exists \text{ sequences } x_k \rightarrow \ox,\; x_k^*\overset{w^*}{\rightarrow} x^* \text{ with } x_k^* \in F(x_k),\;k\in\N\big\}.
\end{eqnarray}

For a {\em proper} extended-real-valued function $\varphi\colon X\rightarrow\overline{\R}:=(-\infty,\infty]$ with $\dom\ph:=\{x\in X \;|\;\varphi(x)<\infty\}\ne\emp$ and for any given $\varepsilon\ge 0$, define the {\em $\varepsilon$-regular subdifferential} of $\varphi$ at $\bar{x}\in \dom \varphi$ by
\begin{equation}\label{FrechetSubdifferential}
\widehat\partial_\varepsilon \varphi(\ox):=\Big\{x^*\in X^*\;\Big|\;\liminf\limits_{x\to \bar x}\frac{\varphi(x)-\varphi(\ox)-\langle x^*,x-\ox\rangle}{\|x-\ox\|}\ge  - \varepsilon\Big\}.
\end{equation}
When $\varepsilon = 0$ in \eqref{FrechetSubdifferential}, this construction is known as the (Fr\'echet) \textit{regular subdifferential}
of $\varphi$ at $\ox$ and is denoted by $\widehat{\partial}\varphi(\ox)$. The (Mordukhovich) {\em limiting/basic subdifferential} of $\varphi$ at $\ox$ is defined via the sequential outer limit \eqref{eq:P-K} by
\begin{equation}\label{MordukhovichSubdifferential}
\partial\varphi(\bar{x}):=\Limsup_{\substack{x\overset{\varphi}{\to}\bar{x},\,\varepsilon \downarrow 0 }}\; \widehat{\partial}_\varepsilon\varphi(x)
\end{equation}
with $x\overset{\varphi}{\to}\bar{x}$ meaning that $x\rightarrow\bar{x}$ and $\varphi(x)\rightarrow\varphi(\bar{x})$. Both constructions $\Hat\partial\varphi(\bar{x})$ and $\partial\varphi(\bar{x})$ reduce to the gradient $\{\nabla\ph(\ox)\}$ if $\ph$ is strictly differentiable at $\ox$ and to the subdifferential of convex analysis 
\begin{equation*}
\partial \varphi (\ox) :=\big\{x^* \in X^*\;\big|\;\la x^*, x-\ox\ra \le \varphi (x) - \varphi (\ox)\;\text{ for all }\;x\in X\big\}
\end{equation*}
when $\ph$ is convex. It has been well recognized that $\ve>0$ can be equivalently dismissed in \eqref{MordukhovichSubdifferential} provided that $\ph$ is lower semicontinuous (l.s.c.) around $\ox$ and the space $X$ is {\em Asplund}, i.e., a Banach space where every separable subspace has a separable dual; this class is sufficiently broad including, in particular, all reflexive spaces. In Asplund spaces, the limiting subdifferential \eqref{MordukhovichSubdifferential} and the related normal and coderivatives constructions for sets and set-valued mappings enjoy comprehensive {\em calculus rules}, which are based on {\em variational/extremal principles} of variational analysis. Some important properties of \eqref{MordukhovichSubdifferential} with using $\ve>0$ in the limiting procedure are also available in general Banach spaces.\vspace*{0.03in}

The next piece of variational analysis needed in what follows concerns {\em variational convergence} of extended-real-valued functions and sets. Denote 
\begin{eqnarray*}
\mathcal{N}_{\infty} &:=&\big\{N\subset\N\;\big|\;\N \setminus N \;\text{ finite}\big\}\;\mbox{ and }\;\mathcal{N}_{\infty}^{\#}:=\big\{N\subset \N\;\big|\;N\;\text{ infinite}\big\}.
\end{eqnarray*} 
Given a sequence of sets $\{C^k\}\subset X$, its {\em outer limit} and {\em inner limit} are defined by  
\begin{eqnarray*}
\Limsup_{k\rightarrow \infty} C^k:&=&\big\{ x\in X\; \big| \; \exists N \in \mathcal{N}^{\#}_{\infty},\;\exists x^k \in C^k\;\text{ with }\;x^k \xrightarrow{N}x\big\},\\
\Liminf_{k\rightarrow \infty} C^k:&=& \big\{ x\in X\; \big| \; \exists N \in \mathcal{N}_{\infty},\ \exists x^k \in C^k\;\mbox{ with }\;x^k \xrightarrow{N}x\big\}.
\end{eqnarray*}
For a sequence $\left\{\ph^k\right\}$ of extended-real-valued functions on $X$, its {\em lower epi-limit} and  {\em upper epi-limit} are defined, respectively, by
$$
\epi (\text{e-}\liminf_{k\rightarrow \infty}\ph^k) := \Limsup_{k\rightarrow \infty} (\epi\ph^k)\;\mbox{ and }\;
\epi (\text{e-}\limsup_{k\rightarrow \infty}\ph^k) := \Liminf_{k\rightarrow \infty} (\epi\ph^k).
$$
If these two functions agree on $X$, then the \textit{epi-limit} function e-$\lim_{k\rightarrow \infty}\ph^k$ exists, i.e., we have
\begin{eqnarray*}
\text{e-}\lim_{k\rightarrow \infty}\ph^k := \text{e-}\liminf_{k\rightarrow \infty}\ph^k = \text{e-}\limsup_{k\rightarrow \infty}\ph^k\;\mbox{ on }\;X.
\end{eqnarray*}
In this case, it is said that the functions $\ph^k$ \textit{epi-converge} to $\ph$, which is denoted by $\ph^k\xrightarrow{e}\ph$. The next characterization of epi-convergence taken from \cite[Proposition~1.14]{Attouch84} plays an important role to verify the characterizations of local maximal monotonicity established in Section~\ref{sec:localsub}.\vspace*{-0.05in}

\begin{Proposition}\label{lem:epi-charac}
Let $X$ be a Banach space, and let $\varphi^k\colon X\to\oR$ as $k\in\N$ and 
$\varphi\colon X\to\oR$ be given functions. Then the following statements are equivalent:

{\bf(i)} $\varphi^k\xrightarrow{e}\varphi$ as $k\to\infty$.

{\bf(ii)} At each point $x\in X$, we have the relationship 
\begin{equation*}
\begin{cases}
\displaystyle\liminf_{k\to\infty}\varphi^k(x^k) \ge \varphi (x) \ \text{ for every sequence }\;x^k \rightarrow x,\\
\disp\limsup_{k\to\infty}\varphi^k (x^k) \le \varphi (x) \ \text{ for some sequence }\;x^k \rightarrow x.
\end{cases}
\end{equation*}
\end{Proposition}\vspace*{-0.05in}

Finally in this section, recall that a function $\varphi:X\rightarrow \overline{\R}$ is {\em weakly sequentially lower semicontinuous} (weakly sequentially l.s.c.) at $\ox \in \dom \varphi$ if for any sequence $\{x_k\}$ which weakly converges to $\ox$, it holds that $\liminf_{k\rightarrow \infty}\varphi (x_k) \ge \varphi (\ox)$. On the other hand, $\varphi$ is {\em weakly lower semicontinuous} (weakly l.s.c.) at $\ox$ if for any $\varepsilon>0$ there exists a neighborhood $U$ of $\ox$ in the weak topology of $X$ such that
\begin{equation*}
\varphi (x)\ge \varphi (\ox)-\varepsilon\;\mbox{ for all }\;x\in U.
\end{equation*}
It is easily to see that the weak l.s.c.\ yields the weak sequential l.s.c., while the reverse implication is not true even in Hilbert spaces; see a counterexample in \cite[Example~2.1]{KYY15}.\vspace*{-0.15in}

\section{Strong Monotonicity and Strong Convexity}\label{str-conv}\vspace*{-0.05in}

This section is manly devoted to the study of various notions of (global) {\em strong monotonicity} for set-valued mappings and {\em strong convexity} for extended-real-valued functions with {\em given moduli}, and to establishing relationships between them in the case of subdifferential operators. 

Recall first that a set-valued mapping $T:X\rightrightarrows X^*$ is (globally) {\em monotone} on $X$ if 
\begin{eqnarray}\label{glob-mon}
\la x_1^* - x_2^*, x_1- x_2 \ra \ge 0 \;\text{ whenever }\; (x_1,x_1^*),(x_2,x_2^*)\in \gph T.
\end{eqnarray}

In what follows, we address the {\em two notions} of strong monotonocity. The first one involves the {\em duality mapping} $J:X\rightrightarrows X^*$ between $X$ and $X^*$ defined by
\begin{equation}\label{duality}
J(x):=\big\{x^*\in X^*\;\big|\;\la x^*,x\ra =\|x\|^2 =\|x^*\|^2\big\},\quad x\in X,
\end{equation}
which plays a significant role in our study. It is well known that $J$ is the subgradient mapping of the quadratic function $\frac{1}{2}\|\cdot\|^2$, and thus it is maximal monotone on $X$ with $\dom J =X$. It is also well known that the duality mapping is {\em single-valued} on $X$ if and only if the space $(X,\|\cdot\|)$ is {\em G\^ateaux smooth}. i.e., its norm is G\^ateaux differentiable on $X\setminus\{0\}$. This subclass of Banach is very large including, after an equivalent renorming, every separable and every reflexive space; see, e.g., the comprehensive book \cite{fabian}, where the reader can find these and other facts from geometric theory of Banach spaces used below. If $X$ is Hilbert, then the duality mapping $J$ reduces to the usual identity mapping $I$ on $X$.

Employing \eqref{duality}, we say that $T\colon X\tto X^*$ on a Banach space $X$ is (globally) {\em strongly monotone} on $X$ with {\em modulus $\sigma>0$} (or {\em $\sigma$-strongly monotone}) if the shift operator $T-\sigma J$ is monotone on $X$, i.e.,
\begin{eqnarray}\label{str-mon}
\la x_1^* - x_2^*, x_1 - x_2 \ra \ge \sigma \la j(x_1)-j(x_2),x_1-x_2\ra 
\end{eqnarray}
whenever $(x_1,x_1^*),(x_2,x_2^*)\in \gph T$ and $(x_1,j(x_1)),(x_2,j(x_2))\in \gph J$. Another notion of strong monotonicity, which is also broadly employed in applications, is defined by using only the norm on $X$. Given $\sigma>0$, we say that $T:X\rightrightarrows X^*$ is {\em norm $\sigma$-strongly monotone} on $X$ if
\begin{equation}\label{norm-mon}
\la x_1^*-x_2^*,x_1-x_2\ra \ge \sigma \|x_1-x_2\|^2\;\mbox{ for all }\;(x_1,x_1^*), (x_2,x_2^*)\in \gph T.
\end{equation}
It is clear that the notions of $\sigma$-strong monotonicity and norm $\sigma$-strong monotonicity agree if $X$ is a Hilbert space, but it is not the case in general. If the modulus $\sigma>0$ is not specified, we simply use the terms ``strong monotonicity" and ``norm strong monotonicity" for the above notions. The {\em maximality} versions of \eqref{glob-mon}, \eqref{str-mon}, and \eqref{norm-mon} are understood as usual with respect to the graph inclusion. Note also the defined monotonicity notions for set-valued mappings induce the corresponding ones for {\em sets} by associating them with the {\em graphs} of such mappings.\vspace*{0.03in}

The following proposition verified in \cite[Proposition~2.2]{KKMP23} is needed in the subsequent material.\vspace*{-0.05in}

\begin{Proposition}\label{prop:preser} Given $T:X\rightrightarrows X^*$, for any $\sigma>0$ we have the assertions:

{\bf(i)} If $T$ is $\sigma$-strongly maximal monotone, then $T-\sigma J$ is maximal monotone on $X$.

{\bf(ii)} If $X$ is reflexive and $T$ is maximal monotone, then $T+\sigma J$ is maximal monotone. 

{\bf(iii)} If $T=\partial \varphi$ for an l.s.c.\ convex function $\varphi:X \rightarrow \overline{\R}$, then $T+\sigma J$ is maximal monotone.

{\bf(iv)} If in the settings of {\rm(ii)} and {\rm(iii)}, the norm $\|\cdot\|$ is G\^ateaux smooth, then the shifted operator $T+\sigma J$ is $\sigma$-strongly maximal monotone.
\end{Proposition}\vspace*{-0.05in}

Next we turn to extended-real-valued functions $\ph\colon X\to\oR$ and consider for them two notions of {\em strong convexity}. For any $\sigma>0$, define {\em quadratically $\sigma$-shifted function} 
\begin{equation}\label{shift}
\psi(x):=\varphi(x) -(\sigma/2)\|x\|^2 \ \text{ for all } \ x\in X
\end{equation}
and say that $\varphi$ is {\em strongly convex} with {\em modulus} $\sigma>0$, or $\sigma$-{\em strongly convex}, if its quadratic shift \eqref{shift} is convex. We also use another notion of strong convexity introduced by Boris Polyak \cite{Polyak} in the form
\begin{equation}\label{ineq-strong}
\varphi\big(\lambda x + (1-\lambda)y\big)+\sigma \lambda (1-\lambda) \dfrac{\|x-y\|^2}{2} \le \lambda \varphi (x) + (1-\lambda)\varphi (y)\;\mbox{for all }\;x,y\in X,\;\lambda\in (0,1),
\end{equation}
which has been studied and applied in many publications; see, e.g., the books \cite{nam,Zalinescu} and the references therein. We label \eqref{ineq-strong} as {\em Polyak $\sigma$-strong convexity} and use the terms ``strong convexity" and ``Polyak strong convexity" for the notions above if the modulus $\sigma>0$ is not specified. 

It is proved in \cite{Nikodem-Pales} that the two strong convexity notions above {\em agree if and only if} $X$ is an {\em inner product space}. It also follows from \cite{Nikodem-Pales} that in any non-Hilbert space, the quadratic function $x\mapsto \|x\|^2$ is $2$-strongly convex but not Polyak 2-strongly convex. On the other hand, the following example in the space $X:=L^1([0,1];\R)$ shows that a function $\ph\colon X\to\oR$ may be Polyak 2-strongly convex while not being 2-strongly convex. Indeed, define $\ph(x):=\int_0^1|x(t)|^2dt$ for $\int_0^1|x(t)|^2dt < \infty$ and $\ph(x):=\infty$ otherwise. To verify that $\ph$ is Polyak 2-strongly convex, consider $u,v\in \dom\ph = L^2([0,1];\R)$ and $\lambda\in(0,1)$. We clearly have for a.e. $t\in[0,1]$ that
\begin{equation}\label{eq:sc-2}
|\lambda u(t)+(1-\lambda)v(t)|^2 +\lambda (1-\lambda)|u(t)-v(t)|^2 \le \lambda|u(t)|^2+(1-\lambda)|v(t)|^2.
\end{equation}
Integrating both sides of \eqref{eq:sc-2} over $[0,1]$ gives us 
\begin{equation*}
\ph\big(\lambda u + (1-\lambda)v\big) + \lambda (1-\lambda)\|u-v\|_{L^2}^2 \le \lambda\ph(u) + (1-\lambda)\ph(v).
\end{equation*}
Therefore, the classical H{\"o}lder inequality ensures that
\begin{equation*}
\ph\big(\lambda u + (1-\lambda)v\big) + \lambda (1-\lambda)\|u-v\|_{X}^2 \le \lambda\ph(u) + (1-\lambda)\ph(v),
\end{equation*}
i.e., $\ph$ is Polyak 2-strongly convex on $X$. To see that the function $\ph$ is not 2-strongly convex on $X$, the reader can check the violation of the underlying inequality $\psi(\lm u+(1-\lm)v)\le\lm\psi(u)+(1-\lm)\psi(v)$ for the quadratic shift \eqref{shift} with $u(t):={\rm log}(t^{100}+0.5)$, $v(t):={\rm log}(1.5-t)$, and $\lm:=0.5$.\vspace*{0.03in}

It follows from \cite[Corollary~3.5.11]{Zalinescu} that for any l.s.c.\ convex function $\ph\colon X\to\oR$ on an arbitrary Banach space and any $\sigma>0$, the Polyak $\sigma$-strong convexity of $\ph$ is {\em equivalent} to the norm $\sigma$-strong monotonicity \eqref{norm-mon} of the {\em subgradient mapping} $\partial\ph$. Regarding relationships between the $\sigma$-strong convexity of $\ph\colon X\to\oR$ defined via the convexity of \eqref{shift} and the $\sigma$-strong monotonicity \eqref{str-mon} of $\partial\ph$ in generality, we present the following proposition, which involves the limiting subgradient mapping \eqref{MordukhovichSubdifferential}.\vspace*{-0.05in}

\begin{Proposition}\label{strong-lem} For a proper extended-real-valued function $\varphi: X\rightarrow\oR$ defined on a Banach space $(X,\|\cdot\|)$ and for any $\sigma>0$, we have the assertions:

{\bf(i)} If $\varphi$ is $\sigma$-strongly convex, then for all $x\in \dom \partial \varphi$, it holds the representation
\begin{equation}\label{repr}
\begin{array}{ll}
\partial\varphi (x)=\big\{x^* \in X^*\;\big|&\exists j(x)\in J(x)\;\mbox{ with }\; \varphi (u) \ge \varphi(x) + \la x^*, u-x\ra\\
&+(\sigma/2)(\|x\|^2 - 2 \la j(x),u\ra + \|u\|^2)\;\mbox{ whenever }\;u\in X\big\}.
\end{array}
\end{equation}

{\bf(ii)} Suppose that $(X,\|\cdot\|)$ is G\^ateaux smooth and that $\varphi$ is l.s.c. If $\varphi$ is $\sigma$-strongly convex, then $\partial \varphi$ is $\sigma$-strongly monotone (in fact, $\sigma$-strongly maximal monotone). The converse implication holds if $X$ is an Asplund space.
\end{Proposition}\vspace*{-0.15in}
\begin{proof} To verify (i), assume that $\varphi$ is $\sigma$-strongly convex, i.e., $\psi$ from \eqref{shift} is convex. Pick $x\in\dom\partial\ph$ and deduce from the Moreau-Rockafellar theorem that
\begin{equation*}
\partial \varphi (x) = \partial \big( \psi + (\sigma/2)\|\cdot\|^2\big)(x) = \partial \psi (x) + \sigma J(x),
\end{equation*}
which ensures the claimed representation in (i) by  the subdifferential construction of convex analysis.

To justify the first assertion in (ii), we clearly have 
\begin{equation*}
\partial \varphi  = \partial\big( \varphi - (\sigma/2)\|\cdot\|^2\big) + \partial\big(\sigma/2)\|\cdot\|^2\big) = \partial\big( \varphi - (\sigma/2)\|\cdot\|^2\big) + \sigma J,
\end{equation*}
and thus deduce the result from Proposition~\ref{prop:preser}(iv). 

It remains to verify the converse statement of (ii) provided that the space $X$ is Asplund. By the $\sigma$-strong monotonicity of $\partial\ph$, we get that the operator $\widehat{\partial}\varphi - \sigma J$ is monotone on $X$. It follows from the difference rule for regular subgradients (see, e.g., \cite[Lemma~6.22]{Mordukhovich18} or \cite[Theorem~3.4]{MNY06}) that
\begin{equation*}
\widehat{\partial}\big(\varphi-(\sigma/2)\|\cdot\|^2\big) \subset \widehat{\partial}\varphi-\sigma J,
\end{equation*}
and hence $\widehat{\partial}\big(\varphi - \frac{\sigma}{2}\|\cdot\|^2\big)$ is monotone. Applying finally \cite[Theorem~3.56]{Mordukhovich06} in Asplund spaces yields the convexity of $\varphi-\frac{\sigma}{2}\|\cdot\|^2$, and thus the function  $\varphi$ is $\sigma$-strongly convex.
\end{proof}\vspace*{-0.15in}

\section{Local Monotonicity of Set-Valued Operators}\label{sec:localmonoproperties}\vspace*{-0.05in}

We begin with defining {\em local monotonicity}, {\em local strong monotonicity}, and their {\em maximality} versions used below. Recall that by definition a ``neighborhood" is an {\em open} set around a given point.\vspace*{-0.05in} 

\begin{Definition}\label{defi:local-mono} Let $T:X\rightrightarrows X^*$, and let $(\ox,\ox^*)\in \gph T$. We say that:

{\bf(i)} $T$ is {\sc locally monotone} around $(\ox,\ox^*)$ if there exist a neighborhood $U\times V$ of $(\ox,\ox^*)$ and a globally monotone operator $\bar{T}:X\rightrightarrows X^*$ such that 
\begin{eqnarray}\label{eq:local-mono}
\gph \bar{T} \cap (U\times V) = \gph T \cap (U\times V).
\end{eqnarray}
$T$ is {\sc locally maximal monotone} around $(\ox,\ox^*)$ if there exist a neighborhood $U\times V$ of $(\ox,\ox^*)$ and a maximal monotone mapping $\bar{T}$ on $X$ such that \eqref{eq:local-mono} holds.
    
{\bf(ii)} $T$ is {\sc locally strongly monotone} around $(\ox,\ox^*)$ with modulus $\sigma >0$ if there are neighborhoods $U$ of $\ox$, $V$ of $\ox^*$ and a globally $\sigma$-strongly monotone operator $\bar{T}:X\rightrightarrows X^*$ such that 
\begin{eqnarray}\label{eq:local-strong}
\gph \bar{T} \cap (U\times V) = \gph T \cap (U\times V).
\end{eqnarray}
$T$ is {\sc locally strongly maximal monotone} around $(\ox,\ox^*)$ with modulus $\sigma >0$ if there exist a neighborhood $U\times V$ of $(\ox,\ox^*)$ and a $\sigma$-strongly maximal monotone mapping $\bar{T}$ such that \eqref{eq:local-strong} holds.
\end{Definition}\vspace*{-0.05in}

We say that $T$ is monotone {\em with respect to} $W\subset X\times X^*$ ($W$ may not be open) if \eqref{eq:local-mono} holds with replacing $U\times V$ by $W$. Similar terminology is used for $\sigma$-strong monotonicity and maximality versions.\vspace*{0.03in}

In our preceding paper \cite{KKMP23}, we established various characterizations of the local maximality monotonicity properties of set-valued operators formulated in Definition~\ref{defi:local-mono}. To present a crucial characterization of local maximal monotonicity, let us recall the needed notions and facts from Banach space theory that can be found, e.g., in \cite{fabian} and the references therein. 

A norm $\|\cdot\|$ on $X$ is {\em Fr\'echet smooth} if it is Fr\'echet differentiable on $X\setminus \{0\}$. This important class of Banach spaces is a bit smaller than their G\^ateaux smooth counterpart, but still includes, under equivalent renorming, all reflexive spaces and spaces with separable duals. A space $(X,\|\cdot\|)$ is {\em strictly convex} if for all $x,y\in S_X:=\{x\in X\mid \|x\|=1\}$, $x\ne y$ and all $\lambda\in(0,1)$, we have
$\|\lambda x+(1-\lambda)y\|<1$. Note that the duality mapping $J$ from \eqref{duality} is single-valued and norm-to-weak$^*$ continuous if $X^*$ is strictly convex. The 2-{\em uniform convexity} of $(X,\||\cdot\|)$ can be defined as
\begin{eqnarray}\label{eq:LWP}
\|x+y\|^2 + c\|x-y\|^2 \le 2\big(\|x\|^2+\|y\|^2\big)\ \text{ for all }\;x,y\in X,
\end{eqnarray}
with some $c\le 1$, where the case of $c=1$ corresponds to Hilbert spaces. It follows from \cite[Corollary~1]{Xu91} that the $2$-uniform convexity of a Banach space is equivalent to the norm strong monotonicity \eqref{norm-mon} of the duality mapping $J$ with some modulus $0<c_1\le 1$ meaning that for every $x,y\in X$, $j(x)\in J(x)$, and $j(y)\in J(y)$ we have the inequality
\begin{equation}\label{LC-strong}
    \la j(x)-j(y),x-y \ra \ge c_1 \|x-y\|^2.
\end{equation}
The reverse inequality to \eqref{eq:LWP} with `$\ge$' replacing `$\le$' and $c\ge 1$ defines 2-{\em uniformly smooth} spaces. The latter class of spaces can be equivalently described as 
\begin{equation}\label{UWP-Lyapunov}
c\|x-y\|^2 \ge \|x\|^2 - 2\la j(x),y\ra + \|y\|^2\;\mbox{ for each }\;x,y\in X\;\mbox{ and }\;j(x)\in J(x).
\end{equation}
If $(X,\|\cdot\|)$ is a G\^ateaux smooth space (i.e., $J(x)=\{j(x)\}$), then
the right-hand side of inequality \eqref{UWP-Lyapunov} is called the {\em Lyapunov functional} and is denoted by $\Lambda(x,y)$. Finally in this brief tour of Banach space geometry, we mention the {\em Kadec-Klee property} of $(X,\|\cdot\|)$ postulating that for any sequence $x_k\st{w}{\to}x$ with $\|x_k\|\to\|x\|$ we have $\|x_k-x\|\to 0$ as $k\to\infty$. 

Having in hand the Banach space notions discussed above, recall a fundamental renorming result saying that any {\em reflexive} Banach space admits an equivalent norm $\|\cdot\|$ such that
\begin{equation}\label{Troyanski-Asplund}
\|\cdot\| \text{ is Fr\'echet smooth, strictly convex, and has the Kadec-Klee property.}
\end{equation}

The following {\em resolvent characterization} of local maximal monotonicity of set-valued operators taken from \cite[Theotem~3.3]{KKMP23} can be viewed as a far-going local counterpart of the classical Minty characterization of global maximal monotonicity. Given a set-valued mapping $T\colon X\tto Y$ between Banach spaces, we say $g\colon X\to Y$ is a {\em single-valued continuous localization} of $T$  around $(\ox,\oy)\in\gph T$ if there exists a neighborhood $U\times V$ of this point such that $\dom g=U$, $g$ is continuous on $U$, and $\gph g=\gph T\cap(U\times V)$.\vspace*{-0.05in}

\begin{Theorem}\label{loc-Minty} Let $T\colon X\tto X^*$ be a set-valued mapping defined on a reflexive Banach space endowed with the equivalent norm $\|\cdot\|$ satisfying the properties listed in \eqref{Troyanski-Asplund}, let $(\ox,\ox^*)\in\gph T$, and let $J\colon X\to X^*$ be the duality mapping \eqref{duality}. Then $T$ is locally maximal monotone around $(\ox,\ox^*)$ if and only if it is locally monotone around $(\ox,\ox^*)$ and its resolvent $(T+\lm J)^{-1}$ admits a single-valued continuous localization around $(\ox^*+\lambda J(\ox),\ox)$ for any $\lm>0$.
\end{Theorem}\vspace*{-0.2in}

\section{Variational Convexity via Subdifferential Graphs} \label{sec:localsub}\vspace*{-0.05in}

In this section, we define the underlying notions of {\em variational} and {\em strong variational convexity} of extended-real-valued functions on Banach spaces by following the pattern of Rockafellar \cite{rtr251} in finite dimensions. The main result of this section provides a characterization of both notions via the corresponding descriptions of the subdifferential graph.\vspace*{-0.05in}

\begin{Definition}\label{def:vsc} Let $\ph\colon X\to\oR$ be an l.s.c.\ function on a Banach space $X$, and let $\ox\in\dom\ph$. Then:

{\bf(i)} $\ph$ is {\sc variationally convex} at $\ox$ for a limiting subgradient $\bar{x}^*\in\partial\varphi(\bar{x})$ if there exist a convex neighborhood $U \times V$ of $(\bar{x},\bar{x}^*)$ and a convex l.s.c.\ function $\widehat{\varphi} \leq \varphi$ on $U$ such that
\begin{equation}\label{eq:vsc1}
\left(U_{\varepsilon} \times V\right)\cap\operatorname{gph}\partial \varphi=(U\times V)\cap\operatorname{gph} \partial \widehat{\varphi} \text{ and } \varphi(x)=\widehat{\varphi}(x) \text{ at the common elements } (x, x^*)
\end{equation}
for some $\e>0$, where $U_{\e}:=\{x\in U\mid \varphi (x) < \varphi (\ox) + \e\}$.

{\bf(ii)} Given $\sigma>0$, we say that $\ph$ is {\sc $\sigma$-strongly  variationally convex} at $\ox$ for $\bar{x}^*\in\partial\varphi(\bar{x})$ if the convex function $\psi$ in {\rm(i)} is assumed to be $\sigma$-strongly convex.

{\bf(iii)} Given $\sigma\ge 0$, we use the term {\sc variational $\sigma$-convexity} for $\ph$ to unify both variational convexity $(\sigma=0)$ and $\sigma$-strong  variational convexity $(\sigma>0)$ cases. 
\end{Definition}\vspace*{-0.05in}

The following auxiliary construction provides a useful technical tool to deal with  strong variational convexity below. For any $\sigma\ne 0$, define the {\em vertical shear mapping} $\Phi_{\sigma}:X\times X^* \rightarrow X\times X^*$ by 
\begin{equation}\label{eq:shear}
\Phi_{\sigma}(x,x^*) = \big(x,x^*+\sigma J(x)\big)\;\text{ whenever }\; (x,x^*)\in X\times X^*,
\end{equation}
which is a homeomorphism from $X\times X^*$ to itself provided that $X$ is Fr\'echet smooth. It is easy to check that for any  $T:X\rightrightarrows X^*$ we have the transformations
\begin{equation}\label{T-T+J}
\gph T \overset{\Phi_{\sigma}}{\longmapsto} \gph (T+\sigma J)\overset{\Phi_{-\sigma}}{\longmapsto}\gph T.
\end{equation}

\begin{Lemma}\label{lem:shift-shear}
Let $X$ be a Fr\'echet smooth space, and let $\varphi:X\rightarrow \overline{\R}$ be a proper function with $\ox\in \dom \varphi$. Given a neighborhood $U\times V$ of $(\ox,\ox^*)\in \gph \partial \varphi$ as well as constants $\e>0$ and $\sigma\neq 0$, there is a convex neighborhood $Q\times W$ of $(\ox,\ox^*-\sigma J(\ox))$ such that we have the implication
\begin{equation}\label{shift:ep}
\big[(x,x^*)\in (Q_{\e/2}\times W)\cap \gph \partial \psi\big] \Longrightarrow \big[(x,x^*+\sigma J(x)) \in (U_{\e}\times V)\cap \gph \partial \varphi\big],
\end{equation}
where $Q_{\e/2}:=\{x\in Q\mid \psi (x) < \psi (\ox) + \e/2\}$ and $U_{\e}$ is taken from Definition~{\rm\ref{def:vsc}}.
\end{Lemma}\vspace*{-0.15in}
\begin{proof}
The quadratic function $\frac{1}{2}\|\cdot\|^2$ is continuously differentiable on $X$ with $\nabla \left( \frac{1}{2}\|\cdot\|^2\right)=J$. The subdifferential sum rule from \cite[Proposition~1.107(ii)]{Mordukhovich06} gives us the equality
\begin{equation}\label{eq:shift1}
\partial \varphi = \partial \psi + \sigma J \ \ \text{ on }\ X,
\end{equation}
which yields $\gph \partial \varphi = \gph (\partial \psi +\sigma J)$. Furthermore, it follows from \eqref{T-T+J} that
\begin{equation}\label{prep6.1}
\Phi_{\sigma}(\gph \partial \psi) = \gph \partial \varphi
\end{equation}
via the homeomorphism $\Phi_{\sigma}$ defined in \eqref{eq:shear}. As $(\ox,\ox^*-\sigma J(\ox)) = \Phi_{-\sigma} (\ox,\ox^*)$, the set $\Phi_{-\sigma}(U\times V)$ is a neighborhood of $(\ox,\ox^*-\sigma J(\ox))$. This allows us to find convex neighborhoods $Q\subset U$ of $\ox$ and $W$ of $\ox^* - \sigma J(\ox)$ such that $Q\times W \subset \Phi_{-\sigma}(U\times V)$ together with $Q\subset B_{r_1}(\ox)$ and $W\subset B_{r_2}(\ox^*-\sigma J(\ox))$, where $r_2>0$ and $r_1>0$ satisfy the estimates
\begin{equation*}
r_1 < \dfrac{\e}{2(\|\ox^*-\sigma J(\ox)\|+r_2)} \quad \text{ and } \quad r_1 (r_1 + 2\|\ox\|)<\dfrac{\e}{\sigma}.
\end{equation*}
It follows from the above that $(\sigma/2)(\|x\|^2 - \|\ox\|^2)<\e/2$ for all $x\in Q$ yielding
\begin{equation*}
\varphi (x) - \varphi (\ox) < \psi (x) - \psi (\ox) + \e/2
\ \text{ for all } \ x\in Q
\end{equation*}
and thus $Q_{\e/2}\subset U_{\e}$. To verify now \eqref{shift:ep}, consider $(x,x^*)\in (Q_{\e/2}\times W)\cap \gph \partial \psi$ and claim that $(x,x^*+\sigma J(x))\in (U_{\e}\times V)\cap \gph \partial \varphi$. Indeed, \eqref{prep6.1} implies that 
\begin{equation}\label{6.1a}
(x,x^*+\sigma J(x)) = \Phi_{\sigma}(x,x^*) \in \Phi_{\sigma}(\gph \partial \psi) = \gph \partial \varphi.
\end{equation}
Since $Q\times W \subset \Phi_{-\sigma}(U\times V)$, we also have 
\begin{equation}\label{6.1b}
(x,x^*+\sigma J(x)) = \Phi_{\sigma}(x,x^*) \in \Phi_{\sigma}(Q\times W) \subset \Phi_{\sigma}(\Phi_{-\sigma}(U\times V)) = U\times V.
\end{equation}
Combining $Q_{\e/2}\subset U_{\e}$ with \eqref{6.1a} and \eqref{6.1b} gives us
\begin{equation*}
(x,x^*+\sigma J(x)) \in (U_{\e}\times V)\cap \gph \partial \varphi,
\end{equation*}
which therefore verifies \eqref{shift:ep} and completes the proof. 
\end{proof}\vspace*{-0.05in}

The next proposition establishes the equivalence between the strong variational convexity of a function and variational convexity of its quadratic shift.\vspace*{-0.05in}

\begin{Proposition}\label{prop:vc-vsc}
Let $\ph\colon X\to\oR$ be an l.s.c.\ function defined on a Fr\'echet smooth space $(X,\|\cdot\|)$, and let $\sigma>0$. Then $\varphi$ is $\sigma$-strongly variationally convex at $\ox$ for $\ox^*\in \partial \varphi (\ox)$ if and only if its $\sigma$-shift
\begin{equation*}
\psi(x): = \varphi(x) -(\sigma/2)\|x\|^2,\quad x\in X,
\end{equation*}
is variationally convex at $\ox$ for $\ox^* - \sigma J(\ox)\in \partial \psi (\ox)$.
\end{Proposition}\vspace*{-0.15in}
\begin{proof}
It follows from \eqref{eq:shift1} that $\ox^* - \sigma J(\ox) \in \partial \psi (\ox)$. The $\sigma$-strong variational  convexity of $\ph$ at $\ox$ for $\ox^*\in \partial \varphi (\ox)$ gives us a neighborhood $U\times V$ of $(\ox,\ox^*)$ and a $\sigma$-strongly convex function $\widehat{\varphi}$ satisfying \eqref{eq:vsc1}. By Lemma~\ref{lem:shift-shear}, we find a convex neighborhood $Q\times W$ of $(\ox,\ox^*-\sigma J(\ox))$ on which \eqref{shift:ep} holds. Considering the convex function $\widehat{\psi}:= \widehat{\varphi}-\frac{\sigma}{2}\|\cdot\|^2$, deduce from the subdifferential sum rule that $\partial \widehat{\varphi} = \partial \widehat{\psi} + \sigma J$. It is easy to see that $\widehat{\psi}$ is convex and l.s.c.\ with $\widehat{\psi}\le \psi$ on $U$. To verify the variational convexity of $\psi$ at $\ox$ for $\ox^*-\sigma J(\ox)$, it suffices to show that
\begin{equation}\label{eq:suff1}
(Q_{\e/2}\times W) \cap \gph \partial \psi = (Q \times W) \cap \gph \partial \widehat{\psi} \text{ and } \psi (x) = \widehat{\psi}(x) \text{ at the common elements } (x,x^*),
\end{equation}
where $Q_{\e/2}$ is taken from Lemma~\ref{lem:shift-shear}. To furnish this, consider $(x,x^*)\in (Q_{\e/2}\times W) \cap \gph \partial \psi$ and deduce from \eqref{shift:ep} the fulfillment of inclusion
\begin{equation}\label{use-lem-1}
(x,x^*+\sigma J(x)) \in (U_{\e}\times V) \cap \gph \partial \varphi.
\end{equation}
It follows from \eqref{use-lem-1} and \eqref{eq:vsc1} that
$$
(x,x^*+\sigma J(x))\in \gph \partial \widehat{\varphi} \quad \text{ and } \quad \varphi (x) = \widehat{\varphi}(x).
$$
The above representation of $\partial\widehat{\varphi} (x)$ implies  that $(x,x^*)\in \gph \partial \widehat{\psi}$ with $\psi (x) = \widehat{\psi}(x)$.

On the other hand, consider $(x,x^*)\in (Q\times W)\cap \gph \partial \widehat{\psi}$. Utilizing the obtained representation of $\partial\Hat\ph$ and the transformation for $\Phi_{\sigma}$ in \eqref{T-T+J}, we have $(x,x^*+\sigma J(x))\in (U\times V)\cap \gph \partial \widehat{\varphi}$. Then \eqref{eq:vsc1} brings us to $(x,x^*+\sigma J(x))\in \gph \partial \varphi$ and $\varphi (x) = \widehat{\varphi}(x)$, and therefore
$$
(x,x^*)\in \gph \partial \psi \quad \text{with} \quad  \psi (x) = \widehat{\psi} (x).
$$
It remains to show that $x\in Q_{\e/2}$, i.e., $\psi (x) < \psi (\ox) +\e/2$. Indeed, it follows from $(x,x^*)\in (Q\times W)\cap \gph \partial \widehat{\psi}$ and the choice of $Q,W$ above that
\begin{equation*}
\begin{array}{ll}
\psi(x) = \widehat{\psi}(x) \le \widehat{\psi}(\ox)-\langle
x^*, \ox - x\rangle   \le \psi (\ox)-\langle
x^*, \ox - x\rangle  
 \le \psi(\ox) + \|x^*\| \cdot \|\ox - x\| \\ 
\le \psi(\ox) + (\|x^*-(\ox^*-\sigma J(\ox))\|+\|\ox^*-\sigma J(\ox)\|)\cdot \|\ox-x\| \le \psi(\ox) + (r_2 + \|\ox^*-\sigma J(\ox)\|) r_1 < \psi(\ox) +\e/2
\end{array}
\end{equation*}
where $r_1,r_2$ are taken from the proof of Lemma \ref{lem:shift-shear}. This yields $(x,x^*)\in (Q_{\e/2}\times W)\cap \gph \partial \psi$ and thus justifies \eqref{eq:suff1}, which justifies the variational convexity of $\psi$ at $\ox$ for $\ox^*-\sigma J(\ox)$. The reverse implication of the proposition can be verified similarly.
\end{proof}\vspace*{-0.05in}

As a part of the proof of the main result of this section, we present the next lemma about the {\em epi-convergence} of extended-real-valued functions.\vspace*{-0.05in}

\begin{Lemma}\label{Lem:fixed1}
Let $X$ be an arbitrary Banach space, and let $\varphi\colon X\to\oR$ and $\varphi^k\colon X\to\oR$ for $k\in\N$ be l.s.c.\ functions. The following assertions hold:

{\bf(i)} If  $\varphi^k \xrightarrow{e}\varphi$ as $k\to\infty$, then
\begin{equation}\label{Att-1}
\inf \varphi \ge \limsup_{k\to\infty} (\inf \varphi^k).   
\end{equation}

{\bf(ii)} If in addition to {\rm(i)}, we have 
$\argmin \varphi \neq \emptyset$ and 
$\argmin \varphi^k \neq \emptyset$ for all $k\in\N$, then
\begin{equation*}
\Limsup_{k\to\infty} (\argmin \varphi^k) \subset \argmin \varphi.
\end{equation*}
Thus $x^k \rightarrow \ox$ as $k\to\infty$ if $\argmin \varphi$ consists of a unique point $\ox$.
\end{Lemma}\vspace*{-0.15in}
\begin{proof} Assertion (i) is proved in \cite[Proposition~2.9]{Attouch84}. To verify assertion (ii), pick $x \in\Limsup_{k\to\infty}(\argmin \varphi^k)$ and find subsequences $\{k_n\}\subset\N$ and $x^{k_n}\in \argmin \varphi^{k_n}$ such that $x^{k_n}\rightarrow x$ as $n\rightarrow \infty$. It follows from $\varphi^k \xrightarrow{e}\varphi$ that $\Limsup_{k\to\infty}(\epi \varphi^k) = \epi\varphi$. Using further
\begin{equation*}
\big(x, \limsup_{n\to\infty}\varphi^{k_n} (x^{k_n})\big) = \lim_{l\rightarrow \infty}(x^{k_{n_l}}, \varphi^{k_{n_l}} (x^{k_{n_l}}))\;\text{ for some subsequence }\;\{k_{n_l}\}\; \text{ of }\;\{k_n\}
\end{equation*}
and $(x^{k_{n_l}}, \varphi^{k_{n_l}} (x^{k_{n_l}}))\in \epi (\varphi^{k_{n_l}})$ as $l\in \N$ tells us that
\begin{equation*}
\big(x, \limsup_{n\to\infty} \varphi^{k_n} (x^{k_n})\big) \in \Limsup_{k\to\infty} (\epi \varphi^k).
\end{equation*}
Therefore, we have $\big(x,\limsup_{n\to\infty} \varphi^{k_n}(x^{k_n})\big)\in \epi \varphi$ meaning that
\begin{equation*}
\varphi (x) \le \limsup_{n\to\infty} \varphi^{k_n} (x^{k_n}).
\end{equation*}
Getting all the above together brings us to 
\begin{equation*}
\varphi (x) \le \limsup_{n\to\infty} \varphi^{k_n} (x^{k_n}) = \limsup_{n\to\infty} (\inf \varphi^{k_n}) \le \limsup_{k\to\infty} (\inf \varphi^k) \le \inf \varphi,
\end{equation*}
where the last inequality follows from \eqref{Att-1}. This yields $x \in \argmin \varphi$ and thus completes the proof.
\end{proof}\vspace*{-0.05in}

Now we get the main characterizations of this sections for both notions of variational convexity and strong variational convexity of functions on Banach spaces. In finite dimensions, the following theorem is contained in \cite[Theorems 1.1 and 1.2]{rtr251} with the proof based on finite-dimensional geometry.\vspace*{-0.05in}

\begin{Theorem}\label{main:VSC}
Let $\varphi:X \to \overline{\R}$ be an l.s.c.\ function defined on a Banach space $X$, and let $\ox^* \in \partial \varphi(\ox)$ be a basic subgradient from \eqref{MordukhovichSubdifferential}. For any $\sigma\ge 0$, consider the following statements:

{\bf(i)} $\varphi$ is variationally $\sigma$-convex at $\bar{x}$ for $\ox^*$. 

{\bf(ii)} There exist a convex neighborhood $U \times V$ of $(\ox, \ox^*)$ and a number $\varepsilon>0$ for which 
\begin{eqnarray}\label{ii}
\begin{array}{ll}
(u, u^*)\in\left(U_{\varepsilon} \times V\right) \cap \gph \partial \varphi \Longrightarrow \exists\,j(u)\in J(u)\;\mbox{ such that for all }\; x\in U\;\mbox{ we have}\\
\varphi(x) \geq \varphi(u)+ \la u^*,x-u\ra +(\sigma/2)\big( \|u\|^2 - 2 \la j(u),x\ra + \|x\|^2\big).
\end{array}
\end{eqnarray}
Then implication {\rm(i)}$\Longrightarrow${\rm(ii)} is valid in general Banach spaces. The converse implication {\rm(ii)}$\Longrightarrow${\rm(i)} fulfills whenever $\sigma\ge 0$ if $X$ is reflexive and if $\varphi$ is weakly sequentially l.s.c.\ around $\ox$.
\end{Theorem}\vspace*{-0.15in}
\begin{proof} We split the proof of the theorem into {\em eight steps} as follows.\\[0.5ex]
{\bf Step~1:} {\em Implication {\rm(i)}$\Longrightarrow${\rm(ii)} holds in an arbitrary Banach space}. Indeed, suppose that $\ph$ is variationally $\sigma$-convex at $\ox$ for $\ox^*$ and find a convex neighborhood $U\times V$ of $(\ox,\ox^*)$ and an
l.s.c.\ function $\widehat{\varphi}\le \varphi$ on $U$, which is convex for $\sigma=0$ and $\sigma$-strongly convex for $\sigma>0$,
such that the conditions in \eqref{eq:vsc1} are satisfied for some $\e>0$. Without loss of generality, we may suppose that $\widehat{\varphi}$ enjoys the above properties on the whole space $X$ (i.e., associate $U$ with $X$), which allows us to employ Proposition~\ref{strong-lem} on $U$. Pick any $(u,u^*)\in (U_{\e}\times V)\cap \gph \partial \varphi$ and get $(u,u^*)\in\gph\partial\widehat{\varphi}$. Applying Proposition~\ref{strong-lem}(i) tells us that
\begin{equation*}
\exists j(u)\in J(u)\;\mbox{ with }\;\widehat{\varphi}(x)\ge \widehat{\varphi}(u) + \la u^*,x-u\ra+ (\sigma/2)\big( \|u\|^2 - 2 \la j(u),x\ra + \|x\|^2 \big)\;\text{ for all }\;x\in U.
\end{equation*}
Since $\widehat{\varphi}\le \varphi$ on $U$ and $\widehat{\varphi}(u)=\varphi(u)$, it follows that 
$$
\varphi(x) \ge \varphi(u) + \la u^*, x - u\ra+ (\sigma/2)\big(\|u\|^2 - 2 \la j(u),x\ra + \|x\|^2\big)
$$ 
whenever $x \in U$, which justifies implication (i)$\Longrightarrow$ (ii).

Next we start verifying the opposite implication (ii)$\Longrightarrow$(i) under the additional assumptions imposed where they are needed. By (ii), take a convex neighborhood $U\times V$ of $(\ox,\ox^*)$ and a number $\e>0$ such that \eqref{ii} holds. If necessary, we shrink $U$ and $V$ to get $U\subset B_{r_1}(\ox)$ and $V\subset B_{r_2}(\ox^*)$ with $r_1 < \varepsilon/(\|\ox^*\|+r_2)$. For any pair $(u,u^*)\in (U_{\e}\times V)\cap \gph \partial \varphi$, choose $\bar{j}(u)\in J(u)$ satisfying \eqref{ii} and define $\widehat{\varphi}:X \rightarrow \overline{\R}$ by 
\begin{eqnarray}\label{repre1}
\widehat{\varphi}(x) := \sup_{(u,u^*)\in (U_{\e}\times V)\cap{\rm\small gph}\,\partial \varphi}\big\{\varphi(u) + \la u^*, x-u \ra+ (\sigma/2)\|u\|^2-\sigma \la \bar{j}(u),x\ra\big\} + (\sigma/2)\|x\|^2.
\end{eqnarray}
It is easy to see that $\widehat{\varphi}$ is l.s.c.\ and $\sigma$-strongly convex with $\widehat{\varphi}(x)>-\infty$
on $X$. Moreover, it follows from \eqref{ii} and \eqref{repre1} that $\widehat{\varphi}\le \varphi$ on $U$.\\[0.5ex]
{\bf Step~2:} {\em For any l.s.c.\ function $\ph\colon X\to\oR$ on a Banach space $X$, we have $(x,x^*)\in (U_{\varepsilon}\times V)\cap \gph \partial \varphi$ if and only if $(x,x^*)\in (U\times V)\cap \gph \partial \widehat{\varphi}$ with $\widehat{\varphi}(x)=\varphi(x)$}.
To verify this, pick $(x,x^*)\in (U_{\varepsilon}\times V)\cap \gph \partial \varphi$ and show first that $\widehat{\varphi}(x)=\varphi(x)$. Indeed, $\varphi(x)$ can be obviously rewritten as 
\begin{equation*}
\varphi(x) = \varphi(x) + \langle x^*, x - x \ra +(\sigma/2)\|x\|^2-\sigma \la \bar{j}(x),x\ra+(\sigma/2)\|x\|^2\;\text{ with }\;(x,x^*)\in (U_{\varepsilon}\times V)\cap \gph \partial \varphi,
\end{equation*}
which yields $\widehat{\varphi}(x) \ge \varphi(x)$. The latter ensures together with $\widehat{\varphi}\leq \varphi$ on $U$ that $\widehat{\varphi}(x)=\varphi(x)$. Further, we deduce from $(x,x^*)\in (U_{\varepsilon}\times V)\cap \gph \partial \varphi$ and $\widehat{\varphi}(x)=\varphi (x)$ that
\begin{eqnarray*}
\widehat{\varphi}(x)+ \langle x^*, x^{\prime} -x\rangle = \varphi(x)+ \langle x^*, x^{\prime} -x\rangle &\le&  \varphi(x)+ \langle x^*, x^{\prime} -x\rangle+(\sigma/2)\left(\|x\|^2-2\la \bar{j}(x),x^\prime\ra + \|x^\prime\|^2\right) \\
&\le& \widehat{\varphi}(x^\prime)\;\mbox{ for any }\;x^{\prime}\in X,
\end{eqnarray*}where the last inequality follows from \eqref{repre1}. Therefore, 
\begin{equation*}
\widehat{\varphi}(x^\prime) \ge \widehat{\varphi}(x) + \langle x^*, x^\prime -x\rangle\;\mbox{ whenever }\;x^\prime \in X,
\end{equation*}
which tells us that $x^* \in \partial \widehat{\varphi}(x)$ by the convexity of $\widehat{\varphi}$.\vspace*{0.03in}

Now we check the opposite implication claimed in Step~2. Take
$(x,x^*)\in(U\times V) \cap \gph \partial  \widehat{\varphi}$
with $\widehat{\varphi}(x)=\varphi(x)$ and get, by the convexity of $\widehat{\varphi}$ and $x^*\in \partial  \widehat{\varphi}(x)$, that
\begin{equation*}
\langle x^*, x^\prime-x\rangle+\widehat{\varphi}(x)\leq \widehat{\varphi}(x^\prime)\;\text{ for all }\;x^\prime\in X.
\end{equation*}
Then we deduce from $\widehat{\varphi}(x)=\varphi (x)$ and $\widehat{\varphi}\le \varphi$ on $U$ that
\begin{equation}\label{1}
\langle x^*, x^\prime - x\rangle  + \varphi(x) \le \varphi(x^\prime) \ \text{ for all } \ x^\prime  \in U,
\end{equation}
which implies in turn that
$$
\liminf_{x^\prime \to x} \frac{\varphi(x^\prime)-\varphi(x)-\langle  x^*,x^\prime -x\rangle}{\|x^\prime  - x\|} \ge 0.
$$
This tells us that $x^*\in \widehat{\partial}\varphi(x)$ and hence yields $x^*\in \partial \varphi(x)$. To finish the proof in Step~2, it remains to check that $\varphi(x) < \varphi(\ox) + \varepsilon$. To verify this, deduce from \eqref{1} and the above choice of $U$ and $V$  that
\begin{equation*}
\begin{array}{ll}
\varphi(x) \le \varphi(\ox)-\langle
x^*, \ox - x\rangle   \le \varphi(\ox) + \|x^*\| \cdot \|\ox - x\| \\ 
\le \varphi(\ox) + (\|x^*-\ox^*\|+\|\ox^*\|)\cdot \|\ox-x\| \le \varphi(\ox) + (r_2 + \|\ox^*\|) r_1 < \varphi(\ox) + \varepsilon,
\end{array}
\end{equation*}
which shows that $x\in U_\e$, and hence we are done in this step.\vspace*{0.03in}

Let $\O\subset U$ be a convex, closed, and bounded set such that $\ox\in{\rm int}\,\O$ and that $\ph$ is weakly sequentially l.s.c.\ on $\O$. Using the indicator function $\delta_\O$ of $\O$, for each $(x,x^*)\in \gph \partial \widehat{\varphi}$ define 
\begin{equation}\label{psi}
\psi_{x,x^*}(y): = \varphi (y) - \widehat{\varphi}(x) - \la x^*, y - x \ra +(1/2)\|y-x\|^2 + \delta_\O (y),\quad y \in X.
\end{equation}
{\bf Step~3:} {\em In any Banach space $X$, we have
$\psi_{x,x^*}(y)\ge 0$, 
{\rm argmin}$\,\psi_{\ox,\ox^*}=\{\ox\}$, and $\min\psi_{\ox,\ox^*}=0$.} To verify these properties, observe that $\la x^*, y-x\ra +\widehat{\varphi}(x) \le \widehat{\varphi}(y) \le \varphi(y)$ for any $y\in U$ and $(x,x^*)\in \gph \partial \widehat{\varphi}$. By $\O\subset U$, we  get that $\psi_{x,x^*}(y) \ge 0$ whenever $y\in X$. The latter ensures that $\psi_{\ox,\ox^*}(y)=0$ if and only if $y=\ox$. This implies in turn that $\argmin\psi_{\ox,\ox^*}=\{\ox\}$ and that $\min \psi_{\ox,\ox^*}=0$.\\[0.5ex]
{\bf Step~4:} {\em Assume that $\ph\colon X\to\oR$ is weakly sequentially l.s.c.\ around $\ox$ and that $X$ is reflexive. Then the function $\psi_{x,x^*}$ from \eqref{psi} attains its absolute minimum on $X$.} It follows from Mazur's weak closure theorem that the set $\O$ under consideration is weakly sequentially closed. Furthermore, we easily deduce from the weak sequential l.s.c.\ of $\ph$ on $\O$ and the structure of $\psi_{x,x^*}$ in \eqref{psi} that the latter function is weakly sequentially l.s.c.\ on the bounded and weakly sequentially closed subset $\O$ of the reflexive Banach space $X$. Then the appropriate version of the Weierstrass existence theorem in this setting ensures that $\psi_{x,x^*}$ attains its absolute minimum on $\O$.\\[0.5ex]
{\bf Step~5:} {\em For a general l.s.c.\ function $\ph\colon X\to\oR$ on a Banach space $X$, the function $\psi_{x,x^*}$ from \eqref{psi} is epicontinuous in $(x,x^*)$ with respect to $\gph \partial\widehat{\varphi}$, i.e., $\psi_{x_k,x_k^*} \xrightarrow{e}\psi_{x,x^*}$ whenever $(x_k,x_k^*)\xrightarrow{{\rm\small gph}\,\partial \widehat{\varphi}}(x,x^*)$ as $k\to\infty$.} To verify this, we use the characterization of epi-convergence presented in  Proposition~\ref{lem:epi-charac}(ii). Take $(x_k,x_k^*)\xrightarrow{{\rm\small gph}\, \partial \widehat{\varphi}}(x,x^*)$ and then fix any $y\in X$ and $y_k\rightarrow y$ as $k\to\infty$. It easily follows from the convexity of the l.s.c.\ function $\widehat{\varphi}$ that $\widehat{\varphi}(x_k)\rightarrow \widehat{\varphi}(x)$ as $k\to\infty$. This yields
\begin{eqnarray*}
\begin{array}{ll}
\disp\liminf_{k\to\infty}\psi_{x_k,x_k^*}(y_k)\ge\disp\liminf_{k\to\infty}
\varphi(y_k)+\disp\liminf_{k\to\infty}\delta_\O(y_k)-\disp\lim_{k\to\infty}\widehat{\varphi}(x_k)\\
+\disp\lim_{k\to\infty}\big[(1/2)\|y_k-x_k\|^2 -\la x^*_k, y_k - x_k\ra \big]\ge\varphi (y) + \delta_\O(y) - \widehat{\varphi}(x) +(1/2)\|y-x\|^2 - \la x^*,y-x\ra=\psi_{x,x^*}(y).
\end{array}
\end{eqnarray*} 
We can also check directly that
\begin{eqnarray*}
\begin{array}{ll}
\disp\limsup_{k\to\infty}\psi_{x_k,x_k^*}(y)=\disp\limsup_{k\to\infty} \varphi(y)+\disp\limsup_{k\to\infty} \delta_\O(y)-\disp\lim_{k\to\infty} \widehat{\varphi}(x_k)\\
+\disp\lim_{k\to\infty}\big[(1/2)\|y-x_k\|^2 -\la x^*_k, y -x_k\ra \big]=\varphi(y) + \delta_\O(y) - \widehat{\varphi}(x) +(1/2)\|y-x\|^2 - \la x^*,y-x\ra =\psi_{x,x^*}(y)
\end{array}
\end{eqnarray*}
and conclude therefore that $\psi_{x_k,x_k^*}\xrightarrow{e}\psi_{x,x^*}$ as $k\to\infty$, which verifies the claimed epicontinuity.\\[0.5ex]
\textbf{Step~6:} {\em Let $\ph$ be a weakly sequentially l.s.c.\  function around $\ox$ defined on a reflexive space $X$, and let the pairs
$(x_k,x_k^*)\in(U\times V)\cap \gph \partial \widehat{\varphi}$ converge to some $(\ox,\ox^*)$ as $k\to\infty$. There for all large $k$, there is a pair $(\tilde{x}_k,\tilde{x}^*_k)\in (U_{\varepsilon}\times V)\cap \gph \partial \varphi$ such that $\widehat{\varphi}(\tilde{x}_k)=\varphi (\tilde{x}_k)$ and $x_k^* - \tilde{x}^*_k \in J(\tilde{x}_k - x_k)$.} Employing the result of Step~5 and remembering from Step~4 that $\argmin\psi_{x_k,x_k^*}\ne\emp$ allow us to deduce from Lemma~\ref{Lem:fixed1}(ii) the existence of $\{\tilde{x}_k\}$ with $\tilde{x}_k \rightarrow \ox$ as $k\to\infty$ and such that
\begin{eqnarray*}
\tilde{x}_k \in \argmin\psi_{x_k,x^*_k} \cap \mathrm{int}\,\O\;\text{ for all large }\;k\in\N. 
\end{eqnarray*}
Furthermore, we get by Lemma~\ref{Lem:fixed1}(i) that
\begin{equation*}
0=\min\psi_{\ox,\ox^*} \ge \limsup_{k\to\infty}(\inf\psi_{x_k,x_k^*}) = \limsup_{k\to\infty} \big(\psi_{x_k,x_k^*}(\tilde{x}_k)\big),
\end{equation*}
which tells us in turn that $\psi_{x_k,x_k^*}(\tilde{x}_k)\dn 0$ as $k\to\infty$ since $\psi_{x_k,x_k^*}(\tilde{x}_k) \ge 0$. Therefore,
\begin{equation*}
\varphi(\tilde{x}_k) - \widehat{\varphi}(x_k)=\psi_{x_k,x_k^*}(\tilde{x}_k) + \la x^*_k, \tilde{x}_k - x_k\ra -(1/2)\|\tilde{x}_k - x_k\|^2<\e/2
\end{equation*}
from which we derive the conditions
\begin{eqnarray*}
\varphi(\tilde{x}_k) &<& \widehat{\varphi}(x_k) +\e/2\le\widehat{\varphi}(\ox) + \e=\varphi(\ox) + \e,
\end{eqnarray*}
and thus $\tilde{x}_k\in U_{\e}$ for large $k$. Since $\tilde{x}_k$ is a minimizer of $\psi_{x_k,x^*_k}$ on $X$, it satisfies $0\in \partial\psi_{x_k,x^*_k}(\tilde x_k)$. Taking into account that $\ph$ is l.s.c.\ while the other functions in \eqref{psi} are locally Lipschitzian around $\ox$, we apply the subdifferential sum rule from \cite[Theorem~2.33]{Mordukhovich06} to $\partial\psi_{x_k,x^*_k}(\tilde x_k)$ giving us
\begin{eqnarray*}
0 \in \partial\psi_{x_k,x^*_k}(\tilde{x}_k) 
\subset \partial \varphi(\tilde{x}_k) - x^*_k + J(\tilde{x}_k-x_k)\;\mbox{ for large }\;k\in\N
\end{eqnarray*}
since $\partial(\|\cdot-x_k\|^2)(\tilde{x}_k) = J(\tilde{x}_k-x_k)$ and $\partial \delta_\O(\tilde{x}_k)=N(\tilde{x}_k,\Omega)=\{0\}$ due to $\ox\in{\rm int}\,\O$. Thus we find $\tilde{x}^*_k \in \partial \varphi(\tilde{x}_k)$ and $j(\tilde{x}_k-x_k) \in J(\tilde{x}_k-x_k)$ with $\tilde{x}^*_k = x^*_k - j(\tilde{x}_k-x_k)$. It follows from $(x_k,x_k^*)\rightarrow (\ox,\ox^*)$ and $\tilde{x}_k \rightarrow \ox$ that $\tilde{x}_k - x_k \rightarrow 0$. By $\|j(\tilde{x}_k-x_k)\| = \|\tilde{x}_k - x_k\|$, this tells us that $j(\tilde{x}_k-x_k) \rightarrow 0$, and hence $\tilde{x}_k^* \rightarrow \ox^*$. Thus we have $(\tilde{x}_k,\tilde{x}_k^*)\in(U_\e\times V)\cap\gph\partial\ph$ for large $k$. Applying finally the statement of Step~2 ensures that $(\tilde{x}_k,\tilde{x}_k^*)\in \gph \partial \widehat{\varphi}$ and $\widehat{\varphi}(\tilde{x}_k)=\varphi(\tilde{x}_k)$ for large $k$, which verifies the claimed assertion.\\[0.5ex]
{\bf Step~7:} {\em We actually have $\tilde{x}_k=x_k$ for large $k\in \N$ in the setting of Step~{\rm 6}.} Indeed, the convexity of the l.s.c.\ function $\widehat{\varphi}$ yields the monotonicity of $\partial\widehat{\varphi}$.
Due to $(x_k,x^*_k),(\tilde{x}_k,\tilde{x}^*_k)\in\gph \partial\widehat{\varphi}$, we get by Step~6 that 
\begin{equation*}
0 \le \la x^*_k - \tilde{x}^*_k , x_k - \tilde{x}_k \ra = -\|x_k-\tilde{x}_k\|^2,
\end{equation*}
which clearly tells us that $\tilde{x}_k=x_k$.\\[0.5ex]
{\bf Step~8:} {\em In the setting of Step~{\rm 7}, there is a neighborhood $(U,\bar{V})$ of $\ox$ and $\ox^*$ such that for all $(x,x^*)\in (U\times \bar{V})\cap \gph \partial \widehat{\varphi}$ it holds $\widehat{\varphi}(x)=\varphi(x)$}. Supposing the contrary  gives us a sequence $\{(x_k,x_k^*)\}\subset\gph\partial \widehat{\varphi}$ converging to $(\ox,\ox^*)$ and
such that $\widehat{\varphi}(x_k)\ne\varphi(x_k)$ for all $k\in\N$. It follows from Steps~6 and 7 that $\widehat{\varphi}(x_k)=\varphi (x_k)$ for large $k$, which is a contradiction that verifies (ii) $\Longrightarrow$(i) and completes the proof.
\end{proof}\vspace*{-0.2in}

\section{Variational Convexity via Subdifferential Local Monotonicity and Moreau Envelopes}\label{sec:vc-Moreau-localmax}\vspace*{0.05in}

In this section, we investigate the notions of variational convexity and its strong counterpart for extended-real-valued functions on Banach spaces from different viewpoints in comparison with Section~\ref{sec:localsub}; namely, by using subdifferential local monotonicity and Moreau envelopes. 

The {\em Moreau envelope} of an l.s.c.\ function $\ph\colon X\to\oR$ on a Banach space $X$ is defined by
\begin{equation}\label{el}
e_\lambda \varphi(x) := \inf_{w\in X}\big\{ \varphi(w) +(2\lambda)^{-1}\|w-x\|^2\big\},\quad x\in X,
\end{equation}   
where $\lm>0$ is a parameter. The corresponding {\em proximal mapping} is given by
\begin{equation}\label{pl}
 P_\lambda \varphi(x) := \text{\rm argmin}_{w\in X}\big\{ \varphi(w) + (2\lambda)^{-1}\|w-x\|^2\big\},\quad x\in X.
\end{equation}

As in \cite{Thibault}, we say that $\ph\colon X\to\oR$ is {\em prox-bounded} if 
\begin{equation}\label{prox-bounded}
\varphi(x)\ge\alpha\|x-\ox\|^2+\beta\;\mbox{ for some }\;\alpha,\beta\in \R\;\mbox{ and }\;\ox\in X.
\end{equation}
This notion can be characterized  in terms of Moreau envelopes.\vspace*{-0.05in}

\begin{Proposition}\label{prop:charac-prox-bound}
For an l.s.c.\ function $\ph\colon X\to\oR$ on a Banach space $X$, the following are equivalent:

{\bf(i)} $\varphi$ is prox-bounded in the sense of \eqref{prox-bounded}.

{\bf(ii)} For some $\lambda >0$ and $x\in X$, we have $e_{\lambda} \varphi (x) >-\infty$.

{\bf(iii)} There is $\lambda_0 >0$ such that $e_{\lambda}\varphi (x)>-\infty$ for any $0<\lambda<\lambda_0$ and any $x\in X$.
\end{Proposition}\vspace*{-0.15in}
\begin{proof} To get (i)$\Longrightarrow$(ii), take $\al,\beta$ in \eqref{prox-bounded} and observe that for any $\lambda>0$ with $\alpha+(2\lambda)^{-1}\ge 0$, we have 
\begin{equation*}
\varphi(y) +(2\lambda)^{-1}\|y-x\|^2 \ge\big(\alpha+(2\lambda)^{-1}\big) \|y-x\|^2 +\beta\;\text{ whenever }\;y\in X.
\end{equation*}
Then (ii) is justified by taking the infimum on both sides of the above inequality with respect to $y\in X$.

Supposing now that (ii) is satisfied, denote 
\begin{equation*}
\lambda_0:= \sup\big\{\lambda >0\;\big|\;\exists x\in X\;\text{ and }\; e_{\lambda}\varphi(x)>-\infty\big\}
\end{equation*}
and get from (ii) that $\lambda_0 >0$. Picking any $0<\lambda<\lambda_0$ and $x\in X$, find
$\bar{\lambda} > \lambda$ and $\ox\in X$ such that $e_{\bar{\lambda}}\varphi(\ox)>-\infty$. Whenever $y\in X$, we have
\begin{eqnarray*}
\begin{array}{ll}
\varphi(y)+(2\lambda)^{-1}\|y-x\|^2 =\varphi(y)+(2\bar{\lambda})^{-1}\|y-\ox\|^2 +(2\lambda)^{-1}\|y-x\|^2 -(2\bar{\lambda})^{-1}\|y-\ox\|^2 \\
\ge e_{\bar{\lambda}}\varphi(\ox) +(2\lambda)^{-1}\|y-x\|^2 -(2\bar{\lambda})^{-1}(\|y-x\|+\|x-\ox\|^2) \\
=e_{\bar{\lambda}}\varphi(\ox) + \big((2\lambda)^{-1}-(2\bar{\lambda})^{-1}\big)\cdot\|y-x\|^2 -\bar{\lambda}^{-1}\|y-x\|\cdot\|x-\ox\|-(2\bar{\lambda})^{-1}\|x-\ox\|^2.
\end{array}
\end{eqnarray*}
By $e_{\bar{\lambda}}\varphi(\ox)>-\infty$ in (ii) and by $\bar{\lambda}>\lambda$, the term $\varphi(y)+(2\lambda)^{-1}\|y-x\|^2$ is bounded from below uniformly in $y\in X$. This tells us that $e_{\lambda}\varphi (x)
>-\infty$ for any $x\in X$ and $\lm\in(0,\lm_0)$, i.e., (iii) holds. The last implication (iii)$\Longrightarrow$ (i) is obvious, and hence the proof is complete.
\end{proof}\vspace*{-0.05in}

To proceed further with functions on general Banach spaces, we need the following modifications of Moreau envelopes and proximal mappings depending on a fixed dual vector $x^*\in X^*$. The {\em $x^*$-tilted Moreau envelope} for $\varphi\colon X\to\oR$ with a parameter $\lm>0$ is defined by
\begin{equation}\label{defi:Motilted}
e^{x^*}_{\lambda} \varphi (x):= \inf_{w\in X}\big\{ \varphi(w) - \la x^*,w\ra +(2\lambda)^{-1}\|w-x\|^2\big\},\quad x\in X.
\end{equation}
It is easy to see that \eqref{defi:Motilted} is the usual Moreau envelope \eqref{el} of the {\em tilted function} $\ph(x)-\la x^*,x\ra$. The corresponding 
{\em $x^*$-tilted proximal mapping} is
\begin{equation}\label{defi:proxtilted}
P_{\lambda}^{x^*}\varphi(x) := \text{\rm argmin}_{w\in X}\big\{ \varphi(w) - \la x^*,w\ra +(2\lambda)^{-1}\|w-x\|^2\big\},\quad x\in X.
\end{equation}
 In the case of Hilbert spaces, we have the relationships (see, e.g., \cite[Lemma~2.2]{Poliquin04}):
\begin{equation}\label{rela:Moreau}
e_{\lambda}^{x^*}\varphi (x) = e_{\lambda}\varphi (x+\lambda x^*) -\langle x^*,x\rangle -(\lambda/2)\|x^*\|^2,\quad
P_{\lambda}^{x^*}\varphi (x) = P_{\lambda}\varphi (x + \lambda x^*).
\end{equation}

The next proposition addresses {\em fixed points} of the  $x^*$-tilted proximal mappings $P_{\lambda}^{x^*}$ in the case where $x^*:=\ox^*$ in \eqref{defi:proxtilted}  is a {\em proximal subgradient} of $\ph$ at $\ox\in\dom\ph$ that is defined as follows: there are positive numbers $r$ and $\varepsilon$ for which
\begin{equation}\label{fixed1}
\varphi(x)\ge \varphi(\ox) + \la \ox^*, x - \ox\ra -(r/2)\|x - \ox\|^2\;\mbox{ for all }\;x \in B_{\e}(\ox).
\end{equation}\vspace*{-0.1in}

\begin{Proposition}\label{prop:Patxbar} Let $\varphi:X\to \oR$ be an l.s.c.\ function on a Banach space with $\ox\in \dom \varphi$, and let $\ox^*\in X^*$. Then we have the equivalent assertions:

{\bf(i)} $\varphi$ is prox-bounded, and $\ox^*$ is a proximal subgradient of $\varphi$ at $\ox$.

{\bf(ii)} $P^{\ox^*}_{\lambda} \varphi (\ox)=\{\ox\}$ for some $\lambda >0$.

{\bf(iii)} $P^{\ox^*}_{\lambda} \varphi (\ox)=\{\ox\}$ for all $\lm>0$ sufficiently small.
\end{Proposition}\vspace*{-0.15in}
\begin{proof} To verify (i)$\Longrightarrow$(ii), fix $\ox^*,r,\ve$ satisfying \eqref{fixed1} and check that we can make \eqref{fixed1} holding for all $x\in X$ if $r$ is properly modified. By the prox-boundedness of $\ph$, it follows from Proposition~\ref{prop:charac-prox-bound}(iii) that $e_{\lambda}\varphi(\ox)>-\infty$ for some $\lambda>0$. Choosing $\bar{\lambda}>0$ to be sufficiently small tells us that
\begin{equation}\label{K3}
\varphi(\ox)+\langle  \ox^*,x -\ox\rangle -(2\bar{\lambda})^{-1}\|x -\ox\|^2 \le  e_\lambda \varphi(\ox) -(2\lambda)^{-1}\|x -\ox\|^2 \le \varphi (x)
\end{equation}
holds for all $x$ with $\|x -\ox\|>\varepsilon$. Denoting $\bar{r}:=2\max \{r,1/\bar{\lambda}\}$, deduce from \eqref{fixed1} and \eqref{K3} that
\begin{equation*}
\varphi(x)-\la \ox^*,x\ra +(\bar{r}/2)\|x-\ox\|^2 > \varphi(\ox) - \la \ox^*,\ox\ra\;\mbox{ whenever }\;x  \neq \ox,
\end{equation*}
which is equivalent to $P_{\bar{\lambda}}^{\ox^*}\varphi(\ox)=\{\ox\}$ with $\bar{\lambda}:=1/\bar{r}$ and hence justifies (ii).

Since implication (ii)$\Longrightarrow$(iii) is obvious, it remains to verify (iii)$\Longrightarrow$(i). Indeed, (iii) ensures by Proposition~\ref{prop:charac-prox-bound}(ii) that $\ph$ is prox-bounded. The fact that $\ox^*$ satisfies \eqref{fixed1} follows from $\ox$ being a global minimizer of the function $\varphi (\cdot) - \la \ox^* , \cdot\ra + (2\lambda)^{-1}\|\cdot - \ox\|^2$ on $X$. 
\end{proof}\vspace*{-0.05in}

Next we recall some important notions of variational analysis broadly used in what follows. They can be found in the books \cite{Rockafellar98} in finite dimensions and \cite{Thibault} in Banach spaces. An l.s.c.\  function $\ph\colon X\to\oR$ is {\em prox-regular} at $\bar{x}\in\dom\ph$ for $\bar{x}^*\in\partial\ph(\ox)$ if there are numbers $\varepsilon>0,r\ge0$ such that
\begin{equation}\label{prox}
\varphi(x)\ge\varphi(u)+\langle u^*,x-u\rangle-(r/2)\|x-u\|^2
\end{equation}
for all $x\in\mathbb{B}_\varepsilon(\bar{x})$ and $(u,u^*)\in\gph\partial\varphi\cap (\mathbb{B}_\varepsilon(\ox)\times\mathbb{B}_\varepsilon(\ox^*))$ with $\varphi(u)<\varphi(\bar{x})+\varepsilon$.
If \eqref{prox} holds for all $x^*\in\partial\varphi(\ox)$, then $\varphi$ is called to be {\em prox-regular at} $\ox$. Further, we say that $\varphi$ is {\em subdifferentially continuous} at $\bar{x}$ for $\ox^*\in \partial\varphi(\bar{x})$ if for any $\varepsilon>0$ there is $\delta>0$ with $|\varphi(x)-\varphi(\bar{x})|<\varepsilon$ whenever $(x,x^*)\in\gph\partial\varphi\cap (\mathbb{B}_\delta(\ox)\times\mathbb{B}_\delta(\ox^*))$. If the latter holds for all $x^*\in\partial\varphi(\ox)$, then $\varphi$ is called to be {\em subdifferentially continuous at} $\ox$. Observe that for subdifferentially continuous functions, the condition ``$\varphi(x)<\varphi(\ox)+\varepsilon$" in the definition of prox-regularity can be dropped. 

The \textit{$\varphi$-attentive $\varepsilon$-localization} of $\partial \varphi$ around $(\ox,\ox^*)\in \gph \partial \varphi$ is the  mapping $T^{\varphi}_{\varepsilon}:X\rightrightarrows X^*$ defined by 
\begin{equation}\label{localization}
\gph T^{\varphi}_{\varepsilon} :=\big\{ (x,x^*)\in \gph \partial \varphi \;\big|\;\|x-\ox\|<\varepsilon,\; |\varphi(x)-\varphi(\ox)|<\varepsilon\;\text{ and }\;\|x^*-\ox^*\|<\varepsilon\big\}.
\end{equation}
If $\varphi$ is l.s.c., a localization can be taken with $\varphi(x) < \varphi(\ox) + \varepsilon$ instead of the absolute value in \eqref{localization}.\vspace*{0.03in}

The proposition below, taken from \cite{Thibault}, concerns the tilted mappings \eqref{defi:Motilted}  and \eqref{defi:proxtilted}.\vspace*{-0.05in}

\begin{Proposition}\label{theo:5.3-Thibault} Let $\varphi: X \rightarrow\oR$ be an l.s.c., prox-bounded, and prox-regular function at $\ox$ for $\ox^*\in \partial \varphi (\ox)$ defined on a $2$-uniformly convex and G\^ateaux smooth space. Then there are $\lambda_0,\varepsilon>0$ such that for any $\lambda\in(0,\lambda_0)$ we can find a neighborhood $U_\lambda$ of $\ox$ on which $e^{\ox^*}_{\lambda}\varphi$ is $\mathcal{C}^1$-smooth with the derivative
\begin{equation}\label{tilt2}
\nabla e_\lambda^{\ox^*} \varphi=\lambda^{-1} J \circ (I-P_\lambda^{\ox^*} \varphi).
\end{equation}
Furthermore, the tilted proximal mapping \eqref{defi:proxtilted} is single-valued and continuous being represented as
$$
P_\lambda^{\ox^*} \varphi(u)=\big(I+\lambda J^{-1} \circ\big(T_{\varepsilon}^\varphi-\ox^*\big)\big)^{-1}(u)\;\text { for all }\;u \in U_\lambda,
$$
where $T_{\varepsilon}^\varphi$ is taken from \eqref{localization} so that $(T_{\varepsilon}^\varphi-\ox^*)(x):=T_{\varepsilon}^\varphi(x)-\ox^*$ is single-valued on $U_\lambda$.
\end{Proposition}\vspace*{-0.02in}

Now we present the needed properties of the duality mapping \eqref{duality}.\vspace*{-0.05in}

\begin{Proposition}\label{lem:J1} Let $(X,\|\cdot\|)$ be a $2$-uniformly convex Banach space whose norm is G\^ateaux smooth. Then the corresponding duality mapping $J$ is single-valued and norm-to-norm continuous on the whole space $X$, i.e., the norm $\|\cdot\|$ is Fr\'echet smooth.
\end{Proposition}\vspace*{-0.15in}
\begin{proof}
Fix $x\in X$ and take a sequence $x_k \to x$ as $k\to\infty$. By the reflexivity of $X$ and the G\^ateaux smoothness of its norm, the dual space $X^*$ is strictly convex, and hence we know from the discussions in Section~\ref{sec:localmonoproperties} that $J$ is single-valued and norm-to-weak$^*$ continuous, which ensures by the reflexivity of $X$ that $J(x_k)\xrightarrow{w} J(x)$ as $k\to\infty$. Moreover, it follows from the convergence $\|x_k\|\to \|x\|$ and definition \eqref{duality} that $\|J(x_k)\|\to\|J(x)\|$ as $k\to\infty$. By taking into account that the imposed assumptions imply the Kadec-Klee property of the space in question, we conclude that $J(x_k)\to J(x)$ as $k\to\infty$, which yields the Fr\'echet smoothness of $(X,\|\cdot\|)$ as discussed above.  
\end{proof}\vspace*{-0.05in}

The following slightly modified notions of monotonicity and strong monotonicity of subgradient operators \cite{rtr251} are used below for characterizing variational and strong variational convexity of functions. Note that the modified notions reduce to those from Definition~\ref{defi:local-mono} for subgradient mappings when the function $\ph$ is subdifferentially continuous at the corresponding point.\vspace*{-0.05in}

\begin{Definition}\label{ph-local mon} Let $\varphi: X \rightarrow\oR$ be an l.s.c.\ function on a Banach space $X$. We say that the subgradient mapping $\partial \varphi: X \rightrightarrows X^*$ is {\sc $\varphi$-locally monotone} $($resp.\ {\sc $\varphi$-locally maximal monotone$)$} around $(\bar{x}, \bar{x}^*)\in \gph \partial\varphi$ if there are a neighborhood $U \times V$ of $(\bar{x}, \bar{x}^*)$ and a number $\ve>0$ such that the mapping $\partial \varphi$ is monotone $($resp. maximal monotone$)$ relative to $U_{\varepsilon} \times V$, where $U_{\varepsilon}$ is taken from Definition~{\rm\ref{def:vsc}}. The notions of {\sc $\varphi$-local strong monotonicity} and its {\sc $\varphi$-local strong maximal monotonicity} counterpart with some modulus $\sigma >0$ are defined similarly.
\end{Definition}\vspace*{-0.05in}

Having in hand the above propositions and definitions, we are now ready to establish the new characterizations of variational and strong variational convexity of extended-real-valued functions on Banach spaces. Our major characterizations are given in the following theorem, which concerns the notion of variational convexity and then will be used for characterizing its strong counterpart. In finite dimensions, the equivalence between assertions (i), (ii), and (iii) of the theorem below has been obtained in \cite[Theorem~1]{rtr251} by a different proof, which is based on finite-dimensional geometry and doesn't require prox-regularity and prox-boundedness of $\ph$. Equivalence (i)$\Longleftrightarrow$(iv) was first established in \cite[Theorem~3.2]{kmp22convex} via the standard Moreau envelope. The last equivalence (i)$\Longleftrightarrow$(v), included in Theorem~\ref{1stequi} due to the proof convenience, is taken from the above Theorem~\ref{main:VSC}, where it is obtained without the prox-regularity and prox-bounded assumptions. \vspace*{-0.05in}

\begin{Theorem} \label{1stequi} Let $\varphi:X\to\oR$ be an l.s.c.\ function around $\ox\in\dom\ph$ defined on a Banach space $X$, and let $\ox^* \in \partial \varphi(\ox)$. Consider the following properties:

{\bf(i)} $\varphi$ is variationally convex at $\ox$ for $\ox^*$.

{\bf(ii)} $\partial \varphi$ is $\varphi$-locally maximal monotone around $(\ox,\ox^*)$.

{\bf(iii)} $\partial \varphi$ is $\varphi$-locally monotone around $(\ox,\ox^*)$.

{\bf(iv)} The $\ox^*$-shifted Moreau envelope $e^{\ox^*}_{\lambda} \varphi$ is locally convex around $\ox$ for all small $\lambda>0$. 

{\bf(v)} There exists a convex neighborhood $U \times V$ of $(\ox, \ox^*)$ along with $\varepsilon>0$ such that
\begin{equation}\label{1-vc}
(u,u^*) \in\left(U_{\varepsilon} \times V\right) \cap \gph \partial \varphi \Longrightarrow \varphi(x) \geq \varphi(u)+ \la u^*,x-u\ra\;\mbox{  for all }\;x\in U.
\end{equation}
Then we have implications {\rm(i)}$\Longrightarrow${\rm(ii)} $\Longrightarrow${\rm(iii)} in general Banach spaces. Implications {\rm(iii)}$\Longrightarrow${\rm(iv)} and {\rm(iv)}$\Longrightarrow${\rm(v)} hold when $X$ is a $2$-uniformly convex space whose norm is G\^ateaux smooth, while 
$\varphi$ is prox-bounded and prox-regular at $\ox$ for $\ox^*$. Lastly, implication {\rm(v)}$\Longrightarrow${\rm(i)} is fulfilled if  the space $X$ is reflexive and the function $\ph$ is weakly sequentially l.s.c.\ around $\ox$.
\end{Theorem}\vspace*{-0.15in}
\begin{proof} 
To verify implication (i)$\Longrightarrow$(ii), we get by Definition~\ref{def:vsc} the existence of a neighborhood $U\times V$ of $(\ox,\ox^*)$ and an l.s.c.\  convex function $\widehat{\varphi}\le\ph$ on $U$ such that \eqref{eq:vsc1} holds with some $\e>0$. The classical Rockafellar theorem  \cite{Rockafellar70-Pac} ensures the global maximal monotonicity of $\partial\widehat{\varphi}$ on $X$. This together with \eqref{eq:vsc1} implies that the mapping $\partial \varphi$ is maximal monotone with respect to $U_{\varepsilon}\times V$, which gives us (ii) by Definition~\ref{ph-local mon}. Implication (ii)$\Longrightarrow$(iii) is obvious.

As a part of verifying implications (iii)$\Longrightarrow$(iv) and (iv)$\Longrightarrow$ (v), observe from Proposition~\ref{theo:5.3-Thibault} that if $X$ is G\^ateaux smooth and $2$-uniformly convex and if $\varphi$ is prox-regular and prox-bounded, then there are positive numbers $\lambda_0,\gamma$ and a $\varphi$-attentive $\gamma$-localization 
\begin{equation}\label{att} 
T_{\gamma}^{\varphi}(x):=\begin{cases}\big\{x^*\in\partial\varphi(x)\big|\;\|x^*-\ox^*\|<\gamma\big\} & \text{if}\quad\|x-\ox\|<\gamma\;\text{ and }\;\varphi(x)<\varphi(\ox)+\gamma,\\
\emptyset  &\text{otherwise}
\end{cases}
\end{equation} 
such that for any $\lambda\in(0,\lambda_0)$ there is a neighborhood $U_\lambda$ of $\ox$ on which $e^{\ox^*}_\lambda\varphi$ is of class $\mathcal{C}^{1}$ with the gradient representation \eqref{tilt2} and
that $P_{\lambda}^{\ox^*}\varphi$ is single-valued and continuous on $U_\lambda$ being represented as 
\begin{equation}\label{tilt1}
P_{\lambda}^{\ox^*}\varphi (u) = \big( I + \lambda J^{-1}\circ (T_{\gamma}^{\varphi}-\ox^*)\big)^{-1}(u)\;\text{ for all }\;u\in U_{\lambda}.
\end{equation}
To verify (iii)$\Longrightarrow$(iv), apply Definition~\ref{ph-local mon} to $\partial\ph$ and find $\varepsilon,r>0$ such that $\partial\ph$ is monotone with respect to the set $W^{\e}:=B_{r}(\ox)^{\e}\times B_r(\ox^*)$, where
$B_r(\ox)^{\varepsilon}:=\{x\in B_r(\ox)\;|\;\varphi(x)<\varphi (\ox) + \varepsilon\}$. Observe that $\ox^*$ is a proximal subgradient of $\varphi$ at $\ox$ by the assumed prox-regularity of $\ph$. Due to the prox-boundedness of $\ph$, we deduce from Proposition~\ref{prop:Patxbar} that
\begin{equation}\label{I-P}
P^{\ox^*}_{\lambda}\varphi(\ox)=\{\ox\}\;\text{ and }\; e_{\lambda}^{\ox^*}\varphi (\ox)=\varphi (\ox)-\langle \ox^*,\ox\rangle.
\end{equation} 
Using the $2$-uniform convexity and  G\^ateaux smoothness of $(X,\|\cdot\|)$ tells us by Proposition~\ref{lem:J1} that $J$ is norm-to-norm continuous. Define the function 
\begin{equation}\label{psi1}
\psi_\lambda(x):=e^{\ox^*}_\lambda\varphi(x) +\big\langle \ox^*,P^{\ox^*}_{\lambda}\varphi (x)\ra-(2\lambda)^{-1}\|P_{\lambda}^{\ox^*}\varphi(x)-x\|^2,\quad x\in U_{\lambda},
\end{equation}
and employ \eqref{tilt2}, \eqref{I-P} to get the equalities
\begin{equation*}
\psi_{\lambda} (\ox) = e^{\ox^*}_\lambda\varphi(\ox) + \langle \ox^*,P^{\ox^*}_{\lambda}\varphi (\ox)\rangle   -(2\lambda)^{-1}\|P_{\lambda}^{\ox^*}\varphi(\ox)-\ox \|^2 = \varphi (\ox),
\end{equation*}
\begin{equation*}
\nabla e_\lambda^{\ox^*} \varphi (\ox)=\lambda^{-1} J \circ (\ox-P_\lambda^{\ox^*} \varphi(\ox))=0.
\end{equation*}
Then we deduce from \eqref{I-P} and the continuity of $\psi_{\lambda}$, $\nabla e_{\lambda}^{\ox^*}\varphi$, and $P^{\ox^*}_{\lambda}\varphi$ around $\ox$ that
\begin{equation}\label{conti1}
P^{\ox^*}_{\lambda}\varphi (x) \in B_{r}(\ox),\quad    \psi_{\lambda}(x) < \varphi (\ox) + \varepsilon, \ \text{ and } \ \|\nabla e_{\lambda}^{\ox^*}\varphi (x)\| < r\;\text{ for all }\;x\in U 
\end{equation}
provided that a neighborhood $U\subset U_{\lambda}$ of $\ox$ is sufficiently small. Using \eqref{defi:Motilted}, \eqref{defi:proxtilted}, and \eqref{psi1} yields
$\psi_\lambda(x)=\varphi (P^{\ox^*}_{\lambda}\varphi (x))$ whenever $x\in U$. Denote $y_i:=P_{\lambda}^{\ox^*}\varphi (x_i)$ for some $x_i\in U$, $i=1,2$, and deduce from \eqref{tilt2}, \eqref{tilt1} and the one-to-one property of the duality mapping that
\begin{equation}\label{tilt22}
\lambda^{-1} J(x_i-y_i) +\ox^* \in T^{\varphi}_{\gamma}(y_i)\;\mbox{ and }\;\nabla e^{\ox^*}_{\lambda}\varphi (x_i)= \lambda^{-1} J(x_i-y_i),\quad i=1,2.
\end{equation}
Then we get by \eqref{conti1} that $y_i\in B_r(\ox)$, $\nabla e_{\lambda}^{\ox^*}\varphi (x_i)+\ox^* \in B_r(\ox^*)$, and
$\varphi (y_i) = \varphi\big(P_{\lambda}^{\ox^*}\varphi (x_i)\big)= \psi_{\lambda} (x_i) <\varphi (\ox) + \varepsilon$ implying that $(y_i,\nabla e^{\ox^*}_{\lambda}\varphi (x_i)+\ox^*)\in W^{\varepsilon} \cap \gph \partial \varphi$ for $i=1,2$. It follows from \eqref{tilt22} that $x_i = y_i+\lambda J^{-1}(\nabla e^{\ox^*}_{\lambda}\varphi (x_i))$, which allows us to deduce from the monotonicity of $\partial \varphi$ with respect to $W^{\varepsilon}$ and the global monotonicity of $J^{-1}$ on $X^*$ that
\begin{eqnarray*}
\begin{array}{ll}
\big\langle \nabla e^{\ox^*}_{\lambda}\varphi (x_1)-\nabla e^{\ox^*}_{\lambda}\varphi (x_2),x_1-x_2\big\rangle =\big\langle  \nabla e^{\ox^*}_{\lambda}\varphi (x_1)-\nabla e^{\ox^*}_{\lambda}\varphi (x_2),y_1-y_2\big\rangle \\\\
+ \lambda\big\langle\nabla e^{\ox^*}_{\lambda}\varphi(x_1)-\nabla e^{\ox^*}_{\lambda}\varphi (x_2),J^{-1}\big(\nabla e^{\ox^*}_{\lambda}\varphi (x_1)\big)-J^{-1}\big(\nabla e^{\ox^*}_{\lambda}\varphi (x_2)\big)\big\rangle\ge 0.
\end{array}
\end{eqnarray*}
The latter ensures that the mapping $\nabla e^{\ox^*}_{\lambda}\varphi$ is locally monotone around $\ox$. Employing now the proof of \cite[Theorem~4.1.4]{Hiriart-Convex} brings us to the local convexity of $e^{\ox^*}_{\lambda}\varphi$ around $\ox$, which justifies (iv).  

Our next goal is to verify (iv)$\Longrightarrow$(v). Suppose
that $X$ is G\^ateaux smooth and $2$-uniformly convex and that $\varphi$ is prox-regular, prox-bounded. Under (iv), we have that $e_{\lambda}^{\ox^*}\varphi$ is convex on $U_\lambda$ for small $\lm>0$. It follows from the proof of \cite[Theorem~4.1.1]{Hiriart-Convex} that
\begin{equation*}
e_{\lambda}^{\ox^*}\varphi(x)\ge e_{\lambda}^{\ox^*}\varphi(u)+\langle\nabla e_{\lambda}^{\ox^*}\varphi(u),x-u\rangle\;\mbox{ for all }\;x,u\in U_\lambda. 
\end{equation*} 
Fix convex neighborhoods $U\subset U_{\lambda}$ of $\ox$ with  $U\subset B_\gg(\ox)$ and $V$ of $\ox^*$ with $V\subset B_\gg(\ox^*)$ such that 
\begin{equation*}
x+\lambda J^{-1}(x^*-\ox^*) \in U_{\lambda}\;\text{ whenever }\;(x,x^*)\in U\times V.
\end{equation*}
Then (v) would follow from the inequality
\begin{equation}\label{convexineq}
\varphi(x^\prime)\geq \varphi(x)+\langle x^*,x^\prime -x \rangle\;\mbox{ for all }\;x^\prime  \in U,\;(x,x^*) \in (U_\gamma\times V)\cap \gph\partial\varphi.
\end{equation}
To verify \eqref{convexineq}, select $x^\prime \in U$, $(x,x^*) \in(U_\gamma\times V)\cap\gph\partial\varphi$ and get from \eqref{att} that $(x,x^*) \in\gph T_{\gamma}^{\varphi}$. Moreover, we deduce from \eqref{tilt2} and \eqref{tilt1} that
\begin{equation*}
x = P^{\ox^*}_{\lambda}\varphi \big( x+\lambda J^{-1}(x^*-\ox^*)\big)\;\mbox{ and }\;\nabla e^{\ox^*}_{\lambda}\varphi\big(x+\lambda J^{-1}(x^*-\ox^*)\big)=x^* -\ox^*,
\end{equation*}
which leads us to the conditions
\begin{eqnarray*}
\begin{array}{ll}
\varphi (x^\prime) - \langle  \ox^*,x^\prime \rangle  + (2\lambda)^{-1}\big\|x'-\big(x'+\lambda J^{-1}(x^*-\ox^*)\big)\big\|^2\ge e^{\ox^*}_{\lambda}\varphi\big(x'+\lambda J^{-1}(x^*-\ox^*)\big) \\
\ge e^{\ox^*}_{\lambda}\varphi\big(x+\lambda J^{-1}(x^*-\ox^*)\big) + \big\langle \nabla e^{\ox^*}_{\lambda}\varphi\big(x+\lambda J^{-1}(x^*-\ox^*)\big), x'-x\big\rangle \\
=\big[\varphi (x) -\langle \ox^*,x\rangle  + (2\lambda)^{-1}\big\|x-\big(x+\lambda J^{-1}(x^*-\ox^*)\big)\big\|^2\big] + \langle x^*-\ox^*,x^\prime  -x\rangle.
\end{array}
\end{eqnarray*}
Thus we arrive at $\varphi(x^\prime)\ge\varphi(x) + \langle x^*, x^\prime  - x\rangle$, which justifies \eqref{convexineq}. The last implication (v)$\Longrightarrow$(i) follows from  Theorem~\ref{main:VSC}, and therefore we are done with the proof of this theorem.
\end{proof}\vspace*{-0.05in}

Next we present a consequence of Theorem~\ref{1stequi} in the case of {\em Hilbert spaces}, where the tilted Moreau envelope can be replaced by the standard one. The following corollary extends to Hilbert spaces the finite-dimensional characterization established in \cite[Theorem~3.2]{kmp22convex}.\vspace*{-0.05in}

\begin{Corollary}\label{coro:VC-Hilbert} 
Let $\varphi: X\to\overline{\R}$ be an l.s.c.\ and prox-bounded function defined on a Hilbert space $X$, and let $\ox^*\in\partial\varphi(\bar{x})$ with $\bar{x}\in\dom\varphi$. Consider the properties:

{\bf(i)} $\varphi$ is variationally convex at $\ox$ for $\ox^*$.

{\bf(ii)} $\varphi$ is prox-regular at $\ox$ for $\ox^*$, and the Moreau envelope $e_\lambda\varphi$ is locally convex around $\bar{x}+\lambda\ox^*$ for all $\lambda>0$ sufficiently small.\\[0.5ex]
Then we always have {\rm(i)}$\Longrightarrow${\rm(ii)}, while {\rm(ii)} $\Longrightarrow${\rm(i)} holds when $\varphi$ is weakly sequentially l.s.c.\ around $\ox$.
\end{Corollary}\vspace*{-0.15in}
\begin{proof} Our intention is to check the equivalence:
\begin{equation}\label{Moreau-shift}
e^{\ox^*}_{\lambda}\varphi\;\text{ is locally convex around }\;\ox \Longleftrightarrow e_{\lambda}\varphi\;\text{ is locally convex around }\;\ox + \lambda\ox^*.
\end{equation}
To verify \eqref{Moreau-shift}, we employ the Hilbert space representation of $e_{\lambda}^{\ox^*}\varphi(x)$ from \eqref{rela:Moreau}. Consider the function $\psi_{\ox^*}(x): = x + \lambda\ox^*$ on $X$ and rewrite \eqref{rela:Moreau} as 
\begin{equation*}
e_{\lambda}^{\ox^*}\varphi(x) = (e_{\lambda}\varphi \circ \psi_{\ox^*})(x)  -\langle \ox^*,x\rangle -(\lambda/2)\|\ox^*\|^2,\quad x\in X,
\end{equation*}
which tells us that the local convexity of $e_{\lambda}\varphi$ around $\ox + \lambda\ox^*$ implies the local convexity of $e^{\ox^*}_{\lambda}\varphi$ around $\ox$. Furthermore, we have from \eqref{rela:Moreau} that
\begin{equation*}
e_{\lambda}\varphi (x) = e^{\ox^*}_{\lambda}\varphi (x-\lambda\ox^*) + \la\ox^*,x\ra -(\lambda/2)\|\ox^*\|^2 =(e^{\ox^*}_{\lambda}\varphi \circ \vt_{\ox^*})(x) + \la\ox^*,x\ra -(\lambda/2)\|\ox^*\|^2
\end{equation*}
via the affine mapping $\vt_{\ox^*}(x):=x-\lambda\ox^*$ on $X$. This ensures that the local convexity of $e^{\ox^*}_{\lambda}\varphi$ around $\ox$ yields this property of $e_{\lambda}\varphi$ around $\ox + \lambda \ox^*$. Thus the claimed result follows from Theorem~\ref{1stequi}.
\end{proof}\vspace*{-0.05in}

Now we proceed with deriving characterizations of {\em  strong variational convexity} of extended-real-valued functions with moduli involved. The following two propositions are instrumental to reduce such characterizations to those for variational convexity established above. \vspace*{0.03in}  

The first proposition provides a relationship between Moreau envelopes of a function and its quadratic shift \eqref{shift}. It is taken from \cite[Lemma~4.3]{kmp22convex}, where the proof holds in Hilbert spaces.\vspace*{-0.05in}

\begin{Proposition}\label{aPhat}
Let $\varphi: X\rightarrow\overline{\R}$ be a prox-bounded l.s.c.\  function, and let $\psi$ be the corresponding $\sigma$-shift with $\sigma\ne 0$. Then for any $\lambda\in(0,|\sigma|^{-1})$ we have 
\begin{equation}\label{Moreau-shift1}
e_{\lambda}\varphi (x) = e_{\lambda/(1+\sigma\lambda)}\psi \left( \dfrac{x}{1+\sigma\lambda}\right)+\dfrac{\sigma}{2(1+\sigma\lambda)}\|x\|^2,\quad x\in X.
\end{equation}
\end{Proposition}\vspace*{-0.05in}

The next proposition shows that the modified local strong  monotonicity of subgradient mappings from Definition~\ref{ph-local mon} reduces to modified ``nonstrong" local monotonicity of their quadratic shifts \eqref{shift}.\vspace*{-0.05in}

\begin{Proposition}\label{prop:shift-mono} Let $\varphi: X\rightarrow\overline{\R}$ be an l.s.c. function defined on a Fr\'echet smooth Banach space $(X,\|\cdot\|)$, and let $\sigma>0$. Then the
subgradient mapping $\partial \varphi$ is $\varphi$-locally strongly monotone around $(\ox,\ox^*)\in \gph \partial\varphi$ with modulus $\sigma$ if and only if the subgradient mapping $\partial\psi$ of the $\sigma$-shifted function $\psi$ from \eqref{shift} is $\psi$-locally monotone around $(\ox,\ox^*- \sigma J(\ox))$.
\end{Proposition}\vspace*{-0.15in}
\begin{proof} To verify the ``only if" part, we find by Definition~\ref{ph-local mon} a neighborhood $U\times V$ of $(\ox,\ox^*)$ and $\e>0$ such that $\partial \varphi$ is $\sigma$-strongly monotone with respect to $U_{\e}\times V$. Then Lemma~\ref{lem:shift-shear} gives us a convex neighborhood $Q\times W$ of $(\ox,\ox^*-\sigma J(\ox))$ on which \eqref{shift:ep} holds. Let us now check that $\partial \psi$ is monotone with respect to $Q_{\e/2}\times W$. Indeed, we get by \eqref{shift:ep} that
\begin{equation*}
(x_i,x_i^*+\sigma J(x_i)) \in (U_{\e}\times V)\cap \gph \partial\varphi\;\mbox{ for all }\;(x_i,x^*_i)\in(Q_{\e/2}\times W)\cap\gph\partial\psi,\;i=1,2.
\end{equation*}
It follows from the $\sigma$-strong monotonicity of $\partial \varphi$ with respect to $U_{\e}\times V$ that 
\begin{equation*}
\la x_1^* + \sigma J(x_1) -\big(x_2^* +\sigma J(x_2)\big), x_1 - x_2 \ra \ge \sigma \la J(x_1)-J(x_2),x_1-x_2\ra,
\end{equation*}
which clearly reduces to the inequality
$\la x_1^* - x_2^*,x_1-x_2\ra \ge 0$ and thus justifies the claimed $\psi$-local monotonicity of $\partial \psi$. The ``if" part is verified similarly.
\end{proof}\vspace*{-0.05in}

Now we are ready to derive the characterizations of strong variational convexity of functions on Banach spaces parallel to those for variational convexity in Theorem~\ref{1stequi}; cf.\ \cite{rtr251} and \cite{kmp22convex} in finite dimensions. 
\vspace*{-0.05in}

\begin{Theorem}\label{2ndequi} Let $\varphi:X \to \overline{\R}$ be an l.s.c.\ function defined on a Banach space $X$, and let $\ox^* \in \partial \varphi(\ox)$ with $\ox\in\dom\ph$. For any $\sigma>0$, consider the assertions:

{\bf(i)} $\varphi$ is $\sigma$-strongly variationally convex at $\ox$ for $\ox^*$.

{\bf(ii)} $\partial \varphi$ is $\varphi$-locally strongly maximal monotone around $(\ox,\ox^*)$ with modulus $\sigma$. 

{\bf(iii)} $\partial \varphi$ is $\varphi$-locally strongly monotone around $(\ox,\ox^*)$ with modulus $\sigma$.

{\bf(iv)} There exists a convex neighborhood $U \times V$ of $(\ox, \ox^*)$ along with $\varepsilon>0$ such that we have \eqref{ii}.\\[0.5ex]
Then {\rm(i)}$\Longrightarrow${\rm(ii)} is valid provided that $(X,\|\cdot\|)$ is G\^ateaux smooth, while {\rm(ii)}$\Longrightarrow$ {\rm(iii)} holds in general. Implication {\rm(iii)}$\Longrightarrow$ {\rm(iv)} is satisfied when $X$ is a $2$-uniformly convex space with a G\^ateaux smooth norm provided that $\ph$ is prox-bounded and the shift
$\psi=\varphi-\frac{\sigma}{2}\|\cdot\|^2$ is 
prox-regular at $\ox$ for $\ox^*-\sigma J(\ox)\in \partial \psi (\ox)$. Finally, {\rm(iv)}$\Longrightarrow${\rm(i)} holds if $X$ is reflexive and $\varphi$ is weakly sequentially l.s.c.\ around $\ox$.
\end{Theorem}\vspace*{-0.15in}
\begin{proof} Assuming (i), we get by Definition~\ref{def:vsc} a neighborhood $U\times V$ of $(\ox,\ox^*)$, an l.s.c.\ $\sigma$-strongly convex function $\widehat{\varphi}:X\rightarrow \overline{\R}$, and a number $\ve>0$ such that 
\begin{equation}\label{vc-1S6}
\gph \partial\widehat{\varphi} \cap (U\times V) = \gph \partial \varphi \cap (U_{\varepsilon}\times V).
\end{equation}
Using Proposition~\ref{strong-lem}(ii) and the G\^ateaux smoothness of $\|\cdot\|$ ensures that $\partial \widehat{\varphi}$ is $\sigma$-strongly maximal monotone. Combining the latter with \eqref{vc-1S6} implies that $\partial \varphi$ is $\sigma$-strongly maximal monotone with respect to $U_{\varepsilon}\times V$, which thus justifies (i)$\Longrightarrow$(ii) due to Definition~\ref{ph-local mon}. Implication (ii)$\Longrightarrow$(iii) is trivial.

Assume now that (iii) holds and that $(X,\|\cdot\|)$ is G\^ateaux smooth and $2$-uniformly convex while $\varphi$ is prox-bounded and $\psi= \varphi - \frac{\sigma}{2}\|\cdot\|^2$ is prox-regular at $\ox$ for $\ox^*-\sigma J(\ox)$. It follows from Proposition~\ref{lem:J1} that the norm $\|\cdot\|$ is Fr\'echet differentiable off the origin. Employing this together with the $\varphi$-local strong monotonicity of $\partial \varphi$ around $(\ox,\ox^*)$ with modulus $\sigma$ in (iii) yields by Proposition~\ref{prop:shift-mono} the $\psi$-local monotonicity of $\partial \psi$ around $(\ox,\ox^*-\sigma J(\ox))$. Since $\psi$ is assumed to be prox-regular at $\ox$ for $\ox^* - \sigma J(\ox)\in \partial \psi (\ox)$ and prox-bounded due to the prox-boundedness of $\ph$, and since $\partial \psi$ is $\psi$-locally monotone around $(\ox,\ox^*-\sigma J(\ox))$, implication (iii) $\Longrightarrow$(v) in Theorem~\ref{1stequi} gives us a convex neighborhood $Q\times W$ of $(\ox, \ox^*-\sigma J(\ox))$ along with a positive number $\delta$ such that
\begin{equation}\label{1st-vc}
\big[(u, u^*) \in (Q^{\psi}_{\delta} \times W) \cap \gph \partial \psi \Longrightarrow \psi (x) \geq \psi (u)+ \la u^*,x-u\ra \big]\;\mbox{ for all }\;x\in Q,
\end{equation} 
where $Q_{\delta}^{\psi}:=\{x\in Q\mid \psi (x) < \psi (\ox) + \delta\}$. By Lemma~\ref{lem:shift-shear}, we find convex neighborhoods $U\subset Q$ of $\ox$ and $V$ of $\ox^*$ ensuring the implication
\begin{equation}\label{shift:ep1}
\big[(x,x^*)\in (U_{\delta/2}\times V)\cap \gph \partial \varphi\big] \Longrightarrow \big[ (x,x^*-\sigma J(x)) \in (Q^{\psi}_{\delta}\times W)\cap \gph \partial \psi\big].
\end{equation}
To verify \eqref{ii}, let $\e:=\delta/2$ and take $(u,u^*)\in (U_{\e}\times V)\cap \gph \partial \varphi$, $x\in U$. Then we get by
\eqref{shift:ep1} that
\begin{equation}\label{Qpsi}
(u,u^*-\sigma J(u)) \in (Q^{\psi}_{\delta}\times W)\cap \gph \partial \psi.
\end{equation}
Since $x\in U\subset Q$, it follows from \eqref{1st-vc} and \eqref{Qpsi} that
\begin{equation*}
\psi (x) \ge \psi (u) + \la u^* - \sigma J(u),x-u\ra,
\end{equation*}
which brings us to \eqref{ii} and hence justifies (iii)$\Longrightarrow$(iv). The last implication (iv)$\Longrightarrow$ (i) is a result of Theorem~\ref{main:VSC}, and thus we are done with the proof.
\end{proof}\vspace*{-0.05in}

It is easy to check that for {\em Hilbert} spaces $X$, the prox-regularity assumption on $\psi$ in Theorem~\ref{2ndequi} is {\em equivalent} to the prox-regularity of $\ph$ at $\ox$ for $\ox^*\in\partial\ph(\ox)$. Furthermore, the Hilbert space setting allows us to establish the following {\em quantitative characterization} of strong variational convexity via the strong convexity of usual Moreau envelope \eqref{el} with modulus interplay.\vspace*{-0.05in}

\begin{Corollary} Let $\varphi:X \to \overline{\R}$ be a prox-bounded function on a Hilbert space $X$ that is weakly sequentially l.s.c.\ around $\ox\in\dom\ph$ and prox-regular at $\ox$ for $\ox^*\in \partial \varphi(\ox)$. Then for any $\sigma>0$, the function $\varphi$ is $\sigma$-strongly variationally convex at $\ox$ for $\ox^*$ if and only if the Moreau envelope $e_{\lambda}\varphi$ is locally $\frac{\sigma}{1+\sigma\lambda}$-strongly convex around $\ox+\lambda \ox^*$ for all $\lm>0$ sufficiently small.
\end{Corollary}\vspace*{-0.15in}
\begin{proof} For the ``only if" part, deduce from Theorem~\ref{2ndequi} and Proposition~\ref{prop:shift-mono} that 
$\partial \psi$ is $\psi$-locally monotone around $(\ox,\ox^*-\sigma \ox)$. By using (iii)$\Longrightarrow$(iv) in Theorem \ref{1stequi}, we know that the tilted Moreau envelope $e_{\theta}^{\ox^* - \sigma \ox}\psi$ is locally convex around $\ox$ for small $\theta>0$. This is equivalent to the local convexity of $e_{\theta}\psi$ around $\ox+\theta (\ox^*-\sigma \ox) = (1-\theta \sigma)\ox + \theta \ox^*$ for such $\theta>0$ due to \eqref{rela:Moreau}. Proposition~\ref{aPhat} yields
\begin{equation}\label{eq:aPhat}
e_{\lambda}\varphi (x)=e_{\lambda/(1+\sigma\lambda)}\psi\left(\dfrac{x}{1+\sigma\lambda}\right)+\dfrac{\frac{\sigma}{1+\sigma\lambda}}{2}\|x\|^2\;\mbox{ for all }\;x\in X\;\mbox{ and small }\;\lambda,
\end{equation}
where $\lambda>0$ can be chosen so that $\lambda/(1+\sigma\lambda)>0$ is also small. Combining this with \eqref{eq:aPhat} and the local convexity of $e_{\theta}\psi$ when $\theta>0$ is sufficiently small, we get that
\begin{equation*}
e_{\lambda}\varphi \text{ is locally }\frac{\sigma}{1+\sigma \lambda}\text{-strongly convex around } \ox + \lambda\ox^*\;\text{ for small }\;\lambda >0.
\end{equation*}
The ``if" part is verified similarly using Proposition \ref{prop:vc-vsc}, Theorem \ref{1stequi} and Proposition \ref{aPhat}.
\end{proof}\vspace*{-0.2in}

\section{Second-Order Characterizations of Variational Convexity}\label{sec:2nd}

In this section, we involve tools of generalized differentiation in variational analysis to provide further characterizations of variational $\sigma$-convexity $(\sigma\ge 0)$ for extended-real-valued functions on {\em Hilbert spaces}, which is assumed unless otherwise stated. Along with the notions of regular and limiting subdifferentials defined in Section~\ref{sec:global-Geometry}, we need some other constructions recalled as follows. 

Let $\O\subset X$ be a nonempty set, and let $\dd(\cdot;\O)$ be its indicator function equal to 0 if $x\in\O$ and $\infty$ otherwise. The {\em regular normal} and {\em limiting normal cone} to $\O$ at $\ox\in\O$ are defined, respectively, as the regular and limiting subdifferentials of the indicator function
\begin{equation*}
\Hat N(\ox;\O):=\Hat\partial\dd(\ox;\O)\;\mbox{ and }\;N(\ox;\O):=\partial\dd(\ox;\O).
\end{equation*}
Note that, while the regular normal cone is always convex but the limiting normal cone is often nonconvex, the latter enjoys much better calculus in both finite and infinite dimensions; see, e.g., \cite{Mordukhovich06,Rockafellar98}. 

Next we recall the two notions of coderivatives for a multifunction $F\colon X\tto Y$ between Hilbert spaces. The {\em regular coderivative}
of $F$ at $(\ox,\oy)\in\gph F$ is defined by
\begin{equation}\label{reg-cod} 
\Hat
D^*F(\ox,\oy)(y^*):=\big\{x^*\in X^*\;\big|\;(x^*,-y^*)\in\Hat N \big((\ox,\oy);\gph F\big)\big\},\quad y^*\in Y 
\end{equation}
via the regular normal cone to the graph of $F$. The {\em mixed coderivative} of $F$ is given in the limiting form
\begin{equation}\label{mix-lim-cod}
D^*_M F(\ox,\oy)(\bar{y}^*):=\Limsup_{(x,y,y^*)\rightarrow (\ox,\oy,\oy^*)}\Hat D^*F(x,y)(y^*),\quad \oy^*\in Y^*,
\end{equation} 
where the outer limit `$\Limsup$' defined in \eqref{eq:P-K} is taken with respect to the norm convergence on $Y^*=Y$ and the weak convergence on $X^*=X$; see \cite{Mordukhovich06} for more details and comprehensive calculus rules available for \eqref{mix-lim-cod} in more general settings. Invoking the coderivatives \eqref{reg-cod} and \eqref{mix-lim-cod} applied to the limiting subdifferential \eqref{MordukhovichSubdifferential}, we get the second-order subdifferential constructions for $\ph\colon X\to\oR$ defined by
\begin{equation}\label{seccombine} 
\breve{\partial}^2\varphi(\bar{x},\ox^*)(u):=(\widehat{D}^*\partial\varphi)(\bar{x},\ox^*)(u)\;\mbox{ and }\;{\partial}^2_M\varphi(\bar{x},\ox^*)(u):=({D}^*_M {\partial}\varphi)(\bar{x},\ox^*)(u),\quad u\in X,
\end{equation} 
known as the {\em combined second-order subdifferential} and the {\em mixed second-order subdifferential} of $\ph$ at $\bar{x}$ relative to $\ox^*\in\partial\ph(\ox)$, respectively. It has been well recognized in variational analysis and optimization that both second-order subdifferentials in \eqref{seccombine} can be explicitly computed for broad classes of extended-real-valued functions, admit reasonable calculus rules (especially for $\partial^2_M\varphi$), provide characterizations of major tilt and full stability properties of optimal solutions and maximal monotonicity of subgradient operators with efficient applications to Newtonian algorithms of nonsmooth optimization and practical models of machine learning, etc.; see, e.g., \cite{CBN,DL,GO16,hos,KKMP23,BorisKhanhPhat,kmptmp,Mordukhovich06,Mordukhovich18,MordukhovichNghia1,mor-roc,mor-sar} and the references therein.\vspace*{0.03in} 

Prior to deriving the main second-order subdifferential characterizations of variational $\sigma$-convexity, we present two technical results of their own interest. The first proposition shows among other things that in the setting under consideration below, the prox-boundedness of the function in question is a consequence of prox-regularity and subdifferential continuity.\vspace*{-0.05in}

\begin{Proposition}\label{prop:prox12} Let $\varphi:X \rightarrow \overline{\R}$ be an l.s.c. function on a Banach space. Assume that $\varphi$ is subdifferentially continuous and prox-regular at $\ox$ for $\ox^*\in\partial\varphi (\ox)$ relative to a neighborhood $U\times V$ of $(\ox,\ox^*)$, and let $\sigma\ge 0$. Denoting $\psi := \varphi  + \delta_U$, we have the following assertions:

{\bf(i)} $\psi$ is prox-bounded on $X$ as well as prox-regular and subdifferentially continuous at $\ox$ for $\ox^*$. 

{\bf(ii)} $\partial \psi (x)=\partial \varphi (x)$ for all $x\in U$. As a consequence, $\gph \partial \psi \cap (U\times V) = \gph \partial \varphi \cap (U\times V)$ and $\gph \partial \psi \cap (U^\psi_{\e}\times V) = \gph \partial \varphi \cap (U^\psi_{\e}\times V)$ for any $\e>0$, where $U^\psi_{\e}:=\{x\in U\mid \psi (x)<\psi(\ox)+\e\}$.

{\bf(iii)} $\breve{\partial}^2 \psi (x,x^*) = \breve{\partial}^2 \varphi (x,x^*)$ and $\partial^2_M \psi (x,x^*) = \partial^2_M \varphi (x,x^*)$ for all $(x,x^*)\in \gph \partial\psi\cap (U\times V)$.

{\bf(iv)} If $X$ is reflexive, $\varphi$ is weakly sequentially l.s.c. around $\ox$, and $\psi$ is $\sigma$-strongly variationally convex at $\ox$ for $\ox^*$, then $\varphi$ is also $\sigma$-strongly variationally convex at $\ox$ for $\ox^*$.
\end{Proposition}\vspace*{-0.15in}
\begin{proof} The prox-regularity and subdiffential continuity of $\psi$ in (i) immediately follow from these properties assumed for $\ph$. Let us verify the prox-boundedness of $\psi$. It follows from the prox-regularity and subdifferential continuity of $\ph$ at $\ox$ for $\ox^*$ and from the construction of $\psi$ that  there is $r\ge 0$ for which
\begin{equation*}
\psi (x) \ge \psi (\ox) + \la \ox^*,x-\ox\ra -(r/2)\|x-\ox\|^2
\end{equation*}
whenever $x\in U$. But since $\psi(x)=\infty$ when $x\in X\setminus U$, the latter holds for all $x\in X$, which yields the prox-boundedness of $\psi$ and hence justifies (i).

The properties in (ii) follow directly from the definitions. To verify (iii), we get from $\gph \partial \psi \cap (U\times V) = \gph \partial \varphi \cap (U\times V)$ in (ii) with the open set $U\times V$ and from \eqref{seccombine} that
\begin{equation*}
\breve{\partial}^2 \psi (x,x^*) =  \breve{\partial}^2 \varphi (x,x^*)\;\mbox{ and }\;\partial^2_M \psi (x,x^*) =  \partial^2_M \varphi (x,x^*)\;\mbox{ for any }\;(x,x^*)\in \gph \partial\psi\cap (U\times V).
\end{equation*}

It remains to verify (iv) under the additional assumptions imposed therein, which ensure characterization \eqref{ii} of variational $\sigma$-convexity in Theorem~\ref{main:VSC}. Applying \eqref{ii} to $\psi$ and remembering the definition of this function, we easily conclude that \eqref{ii} holds for $\ph$ as well, which therefore justifies (iv) by Theorem~\ref{main:VSC} and hence completes the proof of the proposition.
\end{proof}\vspace*{-0.05in}

The next proposition extends the well-known from  finite dimensions to Hilbert spaces.\vspace*{-0.05in}

\begin{Proposition}\label{locally-closed} Let $\varphi: X \rightarrow \overline{\R}$ be an l.s.c. function defined on a Hilbert space $X$, and let $(\ox,\ox^*)\in \gph \partial \varphi$. If $\varphi$ is prox-regular and subdifferentially continuous at $\ox$ for $\ox^*$, then the graphical set $\gph \partial \varphi$ is locally closed around $(\ox,\ox^*)$.
\end{Proposition}\vspace*{-0.15in}
\begin{proof} Proposition~\ref{prop:prox12} allows us to assume that $\varphi$ is prox-bounded. Then it follows from \cite[Proposition~4.4]{Thibault-H} under the subdifferential continuity and prox-regularity of $\varphi$ at $\ox$ for $\ox^*$ that there exist constant $\lambda_0,\e>0$ and an $\e$-localization $T_{\e}$ of $\partial\varphi$ around $(\ox,\ox^*)$ defined by
\begin{equation*}
\gph T_{\e} :=\big\{(x,x^*)\in \gph \partial \varphi\;\big|\; x\in \mathbb{B}_{\e}(\ox)\;\text{ and }\;x^*\in \mathbb{B}_{\e}(\ox^*)\big\}
\end{equation*}
such that for any $0<\lambda < \lambda_0$ we get a closed neighborhood $U_{\lambda}$ of $\ox+\lambda \ox^*$ with
\begin{equation}\label{1'}
P_{\lambda}\varphi = (I+\lambda T_{\e})^{-1}\ \text{ on } U_{\lambda},
\end{equation}
where $P_{\lambda}\varphi$ single-valued and continuous on $U_{\lambda}$. Consider further  the set
\begin{eqnarray*}
W:=\big\{(x,x^*)\in X\times X\;\big|\;x\in \mathbb{B}_{\e}(\ox),\; x^*\in \mathbb{B}_{\e}(\ox^*),\; \text{ and }\;x+\lambda x^* \in U_{\lambda}\}.
\end{eqnarray*}
Clearly, $W=\big(\mathbb{B}_{\e}(\ox) \times \mathbb{B}_{\e}(\ox^*)\big) \cap \vt^{-1} \big(\mathbb{B}_{\e}(\ox)\times U_{\lambda}\big)$ with $\vt (x,x^*):= (x,x+\lambda x^*)$ for all $(x,x^*)\in X\times X$. Since $\vt$ is a homeomorphism on $X\times X$, we get that $W$ is a closed set with $(\ox,\ox^*)\in{\rm int}\,W$. Let us check that $\gph \partial \varphi \cap W$ is closed as well. Indeed, pick $\{(x_k,x^*_k)\}\subset \gph \partial \varphi \cap W$ such that $(x_k,x^*_k) \rightarrow (x,x^*)$ as $k\rightarrow \infty$ and get that $(x,x^*)\in W$. Since $(x_k,x^*_k)\in \gph \partial \varphi \cap W$, we have 
\begin{equation}\label{2}
(x_k,x_k^*) \in \gph T_{\e}\;\mbox{ and }\;x_k+\lambda x^*_k \in U_{\lambda}\;\mbox{ for all }\;k\in\N.
\end{equation}
Combining \eqref{1'} with \eqref{2} and the single-valuedness of $P_{\lambda}\varphi$ on $U_{\lambda}$ yields 
$x_k = P_{\lambda}\varphi (x_k+\lambda x^*_k)$ for all $k\in \N$. Since $P_{\lambda}\varphi$ is continuous on $U_{\lambda}$, we get from the latter that $x = P_{\lambda}\varphi(x+\lambda x^*)$. As $(x,x^*)\in W$, it follows that $x+\lambda x^*\in U_{\lambda}$, and thus from \eqref{1'} we have 
\begin{equation*}
x = P_{\lambda}\varphi (x+\lambda x^*) = (I+\lambda T_{\e})^{-1}(x+\lambda x^*).
\end{equation*}
This immediately implies that $x^*\in T_{\e}(x)$, which yields $(x,x^*)\in \gph \partial\varphi\cap W$ and shows that the set $\gph\partial\varphi$ is locally closed around $(\ox,\ox^*)$. 
\end{proof}\vspace*{-0.05in}

Now we are ready to derive second-order subdifferential characterizations of variational $\sigma$-convexity.\vspace*{-0.05in} 

\begin{Theorem}\label{2ndconvexsubH} Let $X$ be a Hilbert space, and let $\varphi:X\to\overline{\R}$ be subdifferentially continuous at $\bar{x}\in\dom\ph$ for $\ox^*\in\partial\varphi(\bar{x})$ and weakly sequentially l.s.c.\ around $\ox$. Then for any $\sigma\ge 0$, the following are equivalent:

{\bf(i)} $\varphi$ is variationally $\sigma$-convex at $\ox$ for $\ox^*$. 

{\bf(ii)} $\varphi$ is prox-regular at $\ox$ for $\ox^*$ and there exist neighborhoods $U$ of $\bar{x}$ and $V$ of $\ox^*$ such that
\begin{equation}\label{PSDlimit}
\langle z,w\rangle\ge \sigma \|w\|^2\;\text{ whenever}\;z\in\partial^2_M\varphi(x,y)(w),\; (x,y)\in\gph\partial\varphi\cap(U\times V),\;w\in X.
\end{equation} 
				
{\bf(iii)} $\varphi$ is prox-regular at $\ox$ for $\ox^*$ and there exist neighborhoods $U$ of $\bar{x}$, and $V$ of $\ox^*$ such that
\begin{equation}\label{PSDcombinesub}
\langle z,w\rangle\ge \sigma \|w\|^2\;\text{ whenever }\;z\in \breve{\partial}^2\ph(x,y)(w),\;(x,y)\in\gph\partial\varphi\cap(U\times V),\;w\in X.
\end{equation} 
\end{Theorem}\vspace*{-0.15in}
\begin{proof} 
Regarding (i)$\Longrightarrow$(ii), observe first that the prox-regularity of the function $\varphi$ at $\ox$ for $\ox^*$ is a direct consequence of its $\sigma$-strong variational convexity. On the other hand, it follows from Theorem~\ref{1stequi} for $\sigma=0$ (resp.\ Theorem~\ref{2ndequi} for $\sigma>0$) that the variational $\sigma$-convexity and subdifferential continuity of $\varphi$ at $\ox$ for $\ox^*$ yields the local (resp.\ local $\sigma$-strong) maximal monotonicity of $\partial\ph$ around $(\ox,\ox^*)$. Then we get \eqref{PSDlimit} by \cite[Theorem~6.3 and Corollary~6.4]{KKMP23}. 

Implication (ii)$\Longrightarrow$(iii) is trivial since the mixed second-order subdifferential always contains the combined one. To verify the remaining implication (iii)$\Longrightarrow$(i), we deduce from Proposition~\ref{prop:prox12} that (iii) ensures that the function $\psi=\varphi+\delta_U$ is prox-bounded, subdifferentially continuous and prox-regular at $\ox$ for $\ox^*$ while satisfying the condition
\begin{equation}\label{PSDcombinesub-2}
\langle z,w\rangle\ge \sigma \|w\|^2\;\text{ whenever }\;z\in \breve{\partial}^2\psi(x,y)(w),\;(x,y)\in\gph\partial\psi\cap(U\times V),\;w\in X.
\end{equation} 
As follows from Proposition~\ref{locally-closed}, $\gph\partial\psi$ is closed around $(\ox,\ov)$. To apply now the local maximal monotonicity characterizations from \cite[Theorem~6.3~and~Corollary~6.4]{KKMP23} to $\partial\psi$, we need the local hypomonotonicity property of $\partial\psi$ around $(\ox,\ox^*)$, which easily follows from the prox-regularity of $\psi$ at $\ox$ for $\ox^*$ due to Proposition~\ref{prop:prox12}(i). The aforementioned characterizations of the local maximal monotonicity from \cite{KKMP23} tell us that \eqref{PSDcombinesub-2} guarantees the local maximal monotonicity (resp.\ local $\sigma$-strong maximal monotonicity) of $\partial\psi$ around $(\ox,\ox^*)$. Hence we deduce from Theorem~\ref{1stequi} and Theorem~\ref{2ndequi} that $\psi$ is variationally $\sigma$-convex at $\ox$ for $\ox^*$, which together with 
Proposition~\ref{prop:prox12}(iii) justifies the claimed assertion (i).
\end{proof}\vspace*{-0.05in}

Finally, we present the following example, which shows that the {\em subdifferential continuity} is {\em essential} for the characterizations of variational $\sigma$-convexity in Theorem~\ref{2ndconvexsubH} even for {\em univariate} functions.\vspace*{-0.03in}

\begin{Example} {\rm Consider the l.s.c.\ function $\varphi :\R \rightarrow \R$ given by $\varphi(x):=0$ if $x=0$ and $\ph(x):=1$ otherwise. The subgradient mapping $\partial \varphi$ is calculated as 
$\partial \varphi (x)=\R$ for $x=0$ and $\partial \varphi (x)=\{0\}$ otherwise. This shows that $\varphi$ is prox-regular while not subdifferentially continuous at $0$ for $0\in \partial \varphi (0)$. Denoting $U_{1}:= \{x\in \R\mid \varphi (x) < \varphi (0) + 1 \} = \{x\in \R\mid \varphi (x) < 1 \} = \{0\}$, we see that 
\begin{equation*}
(U_1\times \R) \cap \gph \partial \varphi = (\{0\}\times \R)\cap \gph \partial \varphi = \{0\}\times \R,  
\end{equation*}
and thus $\partial\ph$ is $\varphi$-locally strongly monotone around $(0,0)$ with modulus $1$. Theorem~\ref{2ndequi} tells us that $\varphi$ is $1$-strongly variationally convex at $0$ for $0$. Now we check that condition \eqref{PSDcombinesub} fails for $\ph$. Indeed, fix any neighborhood $U\times V$ of $(0,0)$ and pick $(x,0)\in U\times V$ with $x<0$. Locally around $(x,0)$, the graph of $\partial \varphi$ agrees with the graph of the function $f\equiv 0$, and therefore
\begin{equation*}
\breve{\partial}^2\ph(x,0)(1) = (\widehat{D}^*\partial \varphi)(x,0)(1) = (\widehat{D}^* f)(x)(1)=\{0\}.
\end{equation*}
This confirms the failure of \eqref{PSDcombinesub} for $\sigma=1$ and  $w=1$.}
\end{Example}\vspace*{-0.25in}

\section{Strong Variational Convexity and Tilt Stability}\label{sec:tilt}\vspace*{-0.05in}

This section is devoted to establishing relationships between strong variational convexity of extended-real-valued functions on Banach spaces and tilt stability of local minimizers in the associated minimization problems. First we recall the fundamental notion of tilt-stable local minimizers introduced and studied by Poliquin and Rockafellar \cite{Poli} in finite-dimensional spaces.\vspace*{-0.05in}

\begin{Definition}\label{defi:tilt}
For an l.s.c.\ function $\varphi: X \rightarrow \overline{\R}$ on a Banach space, $\bar{x} \in\dom\varphi$ is a {\sc tilt-stable local minimizer} of $\varphi$ if there are neighborhoods $U$ of $\ox$ and $V$ of $0\in X$ such that the mapping
\begin{equation}\label{M-tilt}
M_{V,U}(x^*) :=\argmin\big\{\varphi(x)-\langle x^*, x\rangle\;\big|\;x\in U\big\} \quad \text{for }\ x^*\in V
\end{equation}
is single-valued and Lipschitz continuous with $M_{V,U}(0)=\{\bar{x}\}$. We also say that $\sigma\ge 0$ is a {\sc modulus} of tilt stability for $\ox$ if $\sigma$ is a constant of Lipschitz continuity of \eqref{M-tilt}.
\end{Definition}\vspace*{-0.05in}

Starting with \cite{Poli}, where tilt-stable minimizers of prox-regular and subdifferentially continuous functions $\ph\colon\R^n\to\oR$ were characterized via the positive-definiteness of the {\em basic second-order subdifferential/generalized Hessian} by Mordukhovich $\partial^2\ph(\ox,0):=\partial^2_M\ph(\ox,0)$ in \eqref{seccombine}, so much has been done for characterizations and applications of tilt stability in various classes of optimization problems in finite and infinite dimensions; see \cite{DL,GO16,BorisKhanhPhat,Mordukhovich06,Mordukhovich18,MordukhovichNghia13,mor-roc,mor-sar,rtr251,Thibault} for just a few references. The next theorem establishes close relationships between tilt stability and strong variational convexity in Banach spaces. \vspace*{-0.05in}

\begin{Theorem}\label{vsc-tilt} Let $\varphi:X\rightarrow \overline{\R}$ be an l.s.c. function on a Banach space $(X,\|\cdot\|)$ with $\ox\in \dom \varphi$ and $0\in \partial \varphi (\ox)$. Consider the following assertions:

{\bf(i)} $\varphi$ is strongly variationally convex at $\ox$ for $\ox^*=0$. 

{\bf(ii)} $\ox$ is a tilt-stable local minimizer of $\varphi$.\\[0.5ex]
Then we have implication {\rm(i)}$\Longrightarrow${\rm(ii)} if $X$ is $2$-uniformly convex and G\^ateaux smooth. The reverse implication {\rm(ii)}$\Longrightarrow${\rm(i)} holds if $\varphi$ is weakly sequentially l.s.c.\ around $\ox$, subdifferentially continuous and prox-regular with a constant $r\ge 0$ at $\ox$ for $0$, and if either $X$ is $2$-uniformly smooth with $r<\sigma$ where $\sigma^{-1}$ is the modulus of tilt-stability of $\varphi$ at $\ox$, or $X$ is a Hilbert space.
\end{Theorem}\vspace*{-0.15in}
\begin{proof}
To verify (i)$\Longrightarrow$(ii), suppose that $X$ is $2$-uniformly convex and G\^ateaux smooth. This implies that $\|\cdot\|$ is strictly convex, Fr\'echet smooth by Proposition~\ref{lem:J1}, and Kadec-Klee, i.e., it satisfies \eqref{Troyanski-Asplund}. By definition of strong variational convexity and its characterization in Theorem~\ref{main:VSC}, the strong variational convexity of $\varphi$ at $\ox$ for $\ox^*=0$ yields the existence of $\sigma>0$ and a neighborhood $U\times V$ of $(\ox,0)$ along with a $\sigma$-strongly convex l.s.c.\ function $\widehat{\varphi}:X\rightarrow \overline{\R}$ such that we have \eqref{eq:vsc1} and \eqref{ii} with some $\e>0$.

It follows from Proposition~\ref{strong-lem}(ii) that the $\sigma$-strong convexity of $\widehat{\varphi}$ yields the global $\sigma$-strong maximal monotonicity of the mapping $\partial\widehat{\varphi}$ and hence its $\sigma$-strong maximal monotonicity around $(\ox,0)$. Since $X$ is reflexive and $\|\cdot\|$ satisfies \eqref{Troyanski-Asplund}, deduce from Theorem~\ref{loc-Minty}
that $(\partial \widehat{\varphi})^{-1}$ has a continuous single-valued localization around $(0,\ox)$, i.e., we may shrink the neighborhood $U\times V$ of $(\ox,0)$ such that the mapping $\widehat{T}:V \rightrightarrows U$ defined by $\gph \widehat{T}:= \gph (\partial \widehat{\varphi})^{-1} \cap (V \times U)$ is single-valued and continuous on $V \text{ with } \dom \widehat{T} = V$. The further shrinking of the convex neighborhood $U,V$ if necessary gives us
\begin{equation}\label{chooseU2}
\|u-\ox\| <\sqrt{\e} \ \text{ and }\ \|u^*\| < \sqrt{\e} \ \text{ for all } \ (u,u^*)\in U\times V.
\end{equation}

We split the subsequent proof of (i)$\Longrightarrow$(ii) into the two claims as follows.\\[0.5ex]
{\bf Claim~1:} {\em For $M_{V,U}$ from \eqref{M-tilt}, it follows that}  
\begin{equation}\label{U2}
\gph \partial \varphi \cap (U_{\e}\times V) = \gph M_{V,U}^{-1}.
\end{equation}
Indeed, picking $(u,u^*)\in\gph \partial \varphi \cap (U_{\e}\times V)$ and using \eqref{ii} tell us that
\begin{equation*}
\varphi (x) - \la u^*,x\ra \ge \varphi (u) -\la u^*,u \ra  \ \text{ for all } \ x\in U
\end{equation*}
since the Lyapunov functional $\Lambda(\cdot,\cdot)$ from \eqref{UWP-Lyapunov} always takes nonnegative values. This yields $u\in M_{V,U}(u^*)$ and therefore $(u,u^*)\in \gph M^{-1}_{V,U}$. On the other hand, for any pair $(u,u^*)\in\gph M^{-1}_{V,U}$ we have $(u,u^*)\in U\times V$, and hence the subdifferential Fermat rule from \cite[Proposition~1.114]{Mordukhovich06} implies that
\begin{equation*}
0\in \partial \big(\varphi(\cdot)-\la u^*,\cdot\ra \big)(u) = \partial \varphi (u) - u^*, 
\end{equation*}
i.e., $(u,u^*)\in \gph \partial \varphi$. It follows from $(\ox,0)\in U\times V$ and $u\in M_{V,U}(u^*)$ via definition \eqref{M-tilt} that
\begin{equation*}
\varphi (u) \le \varphi (\ox) + \la u^*, u-\ox\ra < \varphi (\ox) + \e,
\end{equation*}
where the last inequality is due to \eqref{chooseU2}. Hence $(u,u^*)\in\gph \partial \varphi \cap (U_{\e}\times V)$, and we arrive at \eqref{U2}.\\[0.5ex]
{\bf Claim~2:} {\em $M_{V,U}$ is single-valued and Lipschitz continuous on $V$ with $M_{V,U}(0)=\{\ox\}$.} To verify this, we get from $(\ox,0)\in \gph \partial \varphi \cap (U_{\e}\times V)$ and \eqref{U2} that $\ox \in M_{V,U}(0)$. For justifying the single-valuedness of $M_{V,U}$, pick any $x^*\in V$ and let $x:= \widehat{T}(x^*)\in U$. With the aforementioned localization of $(\widehat\partial{\varphi})^{-1}$, we get $(x^*,x)\in \gph (\partial \widehat{\varphi})^{-1}\cap(V\times U)$. Thus it follows from \eqref{eq:vsc1} and \eqref{U2} that $x\in M_{V,U}(x^*)$. Let us now show that the inclusion $(x^*,x_1),(x^*,x_2)\in\gph M_{V,U}$ yields $x_1=x_2$. Indeed, \eqref{U2} tells us that $(x_1,x^*)\in \gph \partial \varphi \cap (U_{\e}\times V)$, and we deduce from \eqref{ii} that 
\begin{equation}\label{nearly}
\varphi (x_2) \ge \varphi (x_1) + \la x^*,x_2-x_1\ra + (\sigma/2)\Lambda (x_1,x_2).
\end{equation}
Since $\varphi (x_1) - \la x^*,x_1\ra = \varphi (x_2)-\la x^*,x_2\ra$ due to $x_1,x_2\in M_{V,U}(x^*)$, it follows from \eqref{nearly} that $\Lambda (x_1,x_2)=0$. The imposed properties of $(X,\|\cdot\|)$ ensures that  $x_1=x_2$. Finally, to check the Lipschitz continuity of $M_{V,U}$ on $V$, take two pairs $(x_1^*,x_1),(x_2^*,x_2)\in\gph M_{{V},U}$ and deduce from \eqref{U2} and \eqref{eq:vsc1} that $(x_i,x_i^*)\in\gph \partial \widehat{\varphi}$ for $i=1,2$. By the global $\sigma$-strong maximal monotonicity of $\partial\widehat{\varphi}$, it follows that
\begin{equation*}
\la x_1^*-x_2^*,x_1-x_2\ra \ge\sigma\la J(x_1)-J(x_2),x_1-x_2\ra.
\end{equation*}
Since $X$ is $2$-uniformly convex, we have \eqref{LC-strong} and therefore
\begin{equation*}
\la x_1^*-x_2^*,x_1-x_2\ra \ge \sigma c_1 \|x_1-x_2\|^2 \ \text{ for some } \ c_1>0.
\end{equation*}
Applying the Cauchy-Schwarz inequality brings us to  
$\|x_1-x_2\|\le (\sigma c_1)^{-1} \|x_1^* - x_2^*\|$, which thus justifies the fulfillment of implication (i)$\Longrightarrow$(ii).

Next we provide the verification of (ii)$\Longrightarrow$ (i). Let $\ox$ be a tilt-stable local minimizer of $\varphi$ with modulus $\sigma^{-1}$, let $\varphi$ be weakly sequentially l.s.c.\ around $\ox$, prox-regular with a constant $r\ge 0$ being subdifferentially continuous at $\ox$ for $0$, and let either $X$ be $2$-uniformly smooth with $r<\sigma$, or $X$ be a Hilbert space. It follows from implication (iv)$\Longrightarrow$ (i) in Theorem~\ref{2ndequi} that for verifying the strong variational convexity of $\varphi$ at $\ox$ for $0$, it suffices to find a convex neighborhood $U\times V$ of $(\ox,0)$ and a number $\bar{\sigma}>0$ such that the implication in \eqref{ii} holds with $\sigma:=\bar\sigma$ for all $x\in U$. Considering first the case where $r<\sigma$ and $X$ is $2$-uniformly smooth, we get the estimate
\begin{equation}\label{UWP-Lyapunov*} 
c\|x-y\|^2 \ge \|x\|^2 - 2\la J(x),y\ra + \|y\|^2\;\text{ for all }\;x,y\in X\;\mbox{ with some }\;c\ge 1.
\end{equation}
Unifying the characterizations of tilt-stable local minimizers from \cite[Theorems~3.2 and 4.2]{MordukhovichNghia13} allows us to find convex neighborhoods $U$ of $\ox$ and $V$ of $0$ ensuring the  {\em second-order growth condition}
\begin{equation*}
\varphi (x) \ge \varphi (u) + \la u^*,x-u\ra +(\sigma/2)\|x-u\|^2\;\mbox{ for all }\;(u,u^*)\in (U\times V) \cap \gph \partial \varphi,\;x\in U.
\end{equation*}
The latter condition together with \eqref{UWP-Lyapunov*} yields the estimate 
\begin{equation*}
\varphi (x) \ge \varphi (u) + \la u^*,x-u\ra +(\sigma/2c)\big(\|u\|^2-2\la J(u),x\ra + \|x\|^2\big),
\end{equation*}
which gives us \eqref{ii} with $\bar{\sigma}:=\sigma/c$ and thus justifies the $\bar{\sigma}$-strong variational convexity of $\varphi$ at $\ox$ for $0$. In the case where $X$ is a Hilbert space, the result follows directly from the tilt stability characterizations in \cite[Theorems~3.2 and 4.2]{MordukhovichNghia13} without using the assumption $r<\sigma$ and involving \eqref{UWP-Lyapunov*}. 
\end{proof}\vspace*{-0.05in}

To conclude this section, we present two remarks about the role of prox-regularity and subdifferential continuity assumptions in the obtained characterizations.\vspace*{-0.05in}

\begin{Remark}\label{rem1} {\rm When $X=\R^n$, it is proved in \cite[Theorem~3]{rtr251}, with the essential usage of finite-dimensional geometry, that the equivalent properties of Theorem~\ref{2ndequi} {\em guarantee} that $\ox$ is a tilt-stable minimizer with modulus $\sigma^{-1}$ of an l.s.c.\ function $\ph$ for $\ox^*=0$ {\em without the prox-regularity and subdifferential continuity} assumptions. However, the {\em reverse} implication may {\em fail} if the prox-regularity of $\ph$ is not imposed. Indeed, this follows from \cite[Example~3.4]{DL} even for the case of univariate functions.}
\end{Remark}\vspace*{-0.1in}

\begin{Remark}\label{rem2}{\rm Due to the equivalence in Theorem~\ref{vsc-tilt}, the second-order subdifferential conditions in (ii) and (iii) of Theorem~\ref{2ndconvexsubH} characterize tilt stability of local minimizers for prox-regular and subdifferentially continuous functions in Hilbert spaces. This has already been known from \cite[Theorem~4.3]{MordukhovichNghia13}. Furthermore, the original characterization of tilt stability for such functions in finite dimensions without the modulus specification was obtained in  \cite[Theorem~3.1]{Poli} in the pointbased form
\begin{equation}\label{tilt-pos}
\la z,w\ra>0\;\mbox{ whenever }\;z\in\partial^2\ph(\ox,0)(w),\;w\ne 0
\end{equation}
via the basic second-order subdifferential $\partial^2\ph$, which is the finite-dimensional version of $\partial^2_M\ph$ from \eqref{seccombine}. It is easy to see that the subdifferential continuity of $\ph$ cannot be avoided in the tilt-stability characterizations of Theorem~\ref{2ndconvexsubH} and of \eqref{tilt-pos} in rather simple situations. Now we show, by some elaboration of \cite[Example~5.4]{DL}, that {\em prox-regularity} of $\ph$ is also {\em essential} for the aforementioned second-order subdifferential characterizations of tilt stability even in the case of {\em continuous} univariate functions.}
\end{Remark}\vspace*{-0.15in}

\begin{Example}{\rm Consider the (nonconvex) function $\ph:[-1,1]\to\R$ defined by 
\begin{equation}\label{phi}
\varphi(x):=
\min\big\{(1+k)^{-1})|x|-(k(k+1))^{-1},k^{-1}\big\}\;\text { for }\;(k+1)^{-1}\leq|x|\le k^{-1}\;\mbox{ for }\;k\in\N
\end{equation}
with $\ph(0):=0$. It is shown in \cite{DL} that $\ox=0$ is a tilt-stable local minimizer of $\varphi$ for any modulus $\kappa >0$, and that $\varphi$ is subdifferentially continuous (in fact, continuous) while not being prox-regular at $\ox$ for $0$. We calculate the limiting subdifferential mapping \eqref{MordukhovichSubdifferential} for \eqref{phi} by
\begin{equation*}
\partial\varphi(x)=\begin{cases}
{\left[0,\frac{k+1}{k}\right]} & \text { if }\;x=\frac{1}{k+1}, \\
{\left[-\frac{k+1}{k}, 0\right]} & \text { if }\;x=-\frac{1}{k+1}, \\
\left\{\frac{k+1}{k}\right\} & \text { if }\;\frac{1}{k+1}<x<\frac{k+2}{(k+1)^2}, \\
\left\{-\frac{k+1}{k}\right\} & \text { if }\;-\frac{k+2}{(k+1)^2}<x<-\frac{1}{k+1}, \\
\left\{0,\frac{k+1}{k}\right\} & \text { if }\;x=\frac{k+2}{(k+1)^2}, \\
\left\{0,-\frac{k+1}{k}\right\} & \text { if }\;x=-\frac{k+2}{(k+1)^2}, \\
\{0\} & \text { if }\;\frac{k+2}{(k+1)^2}<|x|<\frac{1}{k}, \\
{[-1,1]} & \text { if }\;x=0.
\end{cases}
\end{equation*}
Let us show that condition \eqref{PSDcombinesub} fails. For any small $\e>0$, consider the neighborhood $B_{\e}(0)\times B_{\e}(0)$ of $(\ox,0)$. Taking $k\in\N$ with $k^{-1}<\e$ and picking $x\in\big((k+2)/(k+1)^2,k^{-1}\big)$, we have $(x,0)\in \gph \partial \varphi \cap \big(B_{\e}(0)\times B_{\e}(0)\big)$. Defining $r:=\min\{x-(k+2)/(k+1)^2,k^{-1}-x\big\}$ and $f\equiv 0$ gives us
\begin{equation}\label{eq:localization-D*}
\gph \partial \varphi \cap \big( B_{r}(x)\times \R) = \gph f \cap \big( B_{r}(x)\times \R),
\end{equation}
which easily yields the equalities
\begin{equation*}
\breve{\partial}^2 \varphi (x,0)(1) = (\widehat{D}^* \partial \varphi)(x,0)(1)=(\widehat{D}^* f) (x)(1) = \{0\cdot 1\} = \{0\}. 
\end{equation*}
This confirms that condition \eqref{PSDcombinesub} is violated for any small $\e>0$.

Next we show that the pointbased condition \eqref{tilt-pos} also fails for $\ph$ from \eqref{phi}, i.e., 
for any $\kappa>0$ there exists $(w,z)\in \gph \partial^2 \varphi (0,0)$ such that $\la z,w \ra < \kappa\|w\|^2$. To verify this, let $c_k:=0.5\big((k+2)/(k+1)^2+k^{-1}\big)$ for any $k\in\N$. We clearly have that $(c_k,0)\rightarrow (0,0)$ as $k\to\infty$, and that all the points $(c_k,0)$ belong to $\gph \partial \varphi$. Observe that \eqref{eq:localization-D*} also holds with replacing $x$ by $c_k$, and thus 
\begin{equation*}
(\widehat{D}^* \partial \varphi) (c_k,0)(1) = (\widehat{D}^*f)(c_k)(1) = \{0\}\;\mbox{ for all }\;k\in\N.
\end{equation*}
Since $0\in (\widehat{D}^* \partial \varphi) (c_k,0)(1)$ and $(c_k,0)\xrightarrow{{\rm\small gph}\,\partial \varphi} (0,0)$ as $k\to\infty$, we deduce from the definitions that 
$$
0\in(D^* \partial \varphi)(0,0)(1) = \partial^2 \varphi (0,0)(1). 
$$
As $\la 0,1 \ra < \kappa\|1\|^2$ for all $\kappa >0$, condition \eqref{tilt-pos} fails for the function $\ph$ under consideration.}
\end{Example}\vspace*{-0.2in}

\section{Concluding Remarks}\label{sec:conc}\vspace*{-0.05in}

This paper provides a systematic study of variational and strong variational convexity of nonsmooth functions in Banach spaces. Although these notions have been recently introduced and studied in finite dimensions, their importance has been already realized in variational analysis, optimization, and applications. We derive various characterizations of these and related notions in infinite-dimensional spaces and establish their relationships between local maximal monotonicity of subgradient operators, positive-semidefiniteness of second-order subdifferentials, tilt stability of local minimizers, etc. 

Among major open questions of our future research, we mention applications of variational and strong variational convexity to {\em variational sufficiency} in structured optimization, which has been recently studied in finite dimensions \cite{kmp22convex,rtr251,roc,r22,ding}. Then we aim at applying the obtained results to numerical algorithms of nonsmooth optimization in finite-dimensional and infinite-dimensional spaces. 


\vspace*{-0.1in}

\small \begin{thebibliography}{10} 

\bibitem{Attouch84} H. Attouch, {\em Variational Convergence for Functions and Operators}, Pitman, Boston, 1984.

\bibitem{Bauschke2011} H. H. Bauschke and  P. L. Combettes, {\em Convex Analysis and Monotone Operator Theory in Hilbert Spaces}, 2nd edition. Springer, New York, 2017.

\bibitem{Thibault-H} F. Bernard and L. Thibault, Prox-regular functions in Hilbert spaces, {\em J. Math. Anal. Appl.} {\bf 303} (2005), 1--14.

\bibitem{CBN} N. H. Chieu, G. M. Lee, B. S. Mordukhovich and T. T. A. Nghia, Coderivative characterizations of maximal monotonicity for set-valued mappings, {\em J. Convex Anal.} {\bf 23} (2016), 461--480.

\bibitem{DL}  D. Drusvyatskiy and A. S. Lewis, Tilt stability, uniform quadratic growth, and strong metric regularity of the subdifferential, {\em SIAM J. Optim.} {\bf 23} (2013), 256--267.

\bibitem{fabian} M. Fabian, P. Habala, P. H\'ajeck, V. Montesinos S., J. Pelant and V. Zizler, {\em Functional Analysis and Infinite-Dimensional Geometry}, 2nd edition, Springer, New York, 2011.

\bibitem{GO16} H. Gfrerer and J. V. Outrata, On computation of limiting coderivatives of the normal-cone mapping to inequality systems and their applications, {\em Optimization} {\bf 65} (2016), 671--700.

\bibitem{Poliquin04} W. L. Hare and R. A. Poliquin, Prox-regularity and stability of the proximal mapping, {\em J. Convex Anal.} {\bf 14} (2007), 589--606.

\bibitem{hos} R. Henrion, J. V. Outrata and T. Surowiec, Analysis of M-stationary points to an EPEC modeling oligopolistic competition in an electricity spot market, {\em ESAIM:COCV Control} {\bf 18} (2012), 295--317. 

\bibitem{Hiriart-Convex} J.-B. Hiriart-Urruty and C. Lemar\'echal, {\em Fundamentals of Convex Analysis}, Springer, Berlin, 1993.

\bibitem{KKMP23} P. D. Khanh, V. V. H. Khoa, B. S. Mordukhovich and V. T. Phat, Local maximal monotonicity in variational analysis and optimization, arXiv:2308.14193 (2023).

\bibitem{BorisKhanhPhat} P. D. Khanh, B. S. Mordukhovich and V. T. Phat, A generalized Newton method for subgradient systems, {\em Math. Oper. Res.}, DOI: 10.1287/moor.2022.1320 (2022).

\bibitem{kmp22convex} P. D. Khanh, B. S. Mordukhovich and V. T. Phat, Variational convexity of functions and variational sufficiency in optimization, {\em SIAM J. Optim.} {\bf 33} (2023), 1121--1158.
			
\bibitem{kmptmp} P. D. Khanh, B. S. Mordukhovich, V. T. Phat and D. B. Tran, Globally convergent coderivative-based generalized Newton methods in nonsmooth optimization, {\em Math. Program.}, DOI: 10.1007/s10107-023-01980-2 (2023).

\bibitem{KYY15} P. D. Khanh, J.-C. Yao and N. D. Yen, The Mordukhovich subdifferentials and directions of descent, {\em J. Optim. Theory Appl.} {\bf 172} (2017), 518--534.

\bibitem{MNY06} B. S. Mordukhovich, N. M. Nam and N. D. Yen, Fr\'echet subdifferential calculus and optimality conditions in nondifferentiable programming, {\em Optimization} {\bf 55} (2006), 685--708.

\bibitem{Mordukhovich06} B. S. Mordukhovich, {\em Variational Analysis and Generalized Differentiation, I: Basic Theory}, Springer, Berlin, 2006.

\bibitem{Mordukhovich18} B. S. Mordukhovich, {\em Variational Analysis and Applications}, Springer, Cham, Switzerland, 2018.

\bibitem{nam} B. S. Mordukhovich and N. M. Nam, {\em Convex Analysis and Beyond, I: Basic Theory}, Springer, Cham, Switzerland, 2022.

\bibitem{MordukhovichNghia13} B. S. Mordukhovich and T. T. A. Nghia, Second-order variational analysis and characterizations of tilt-stable optimal solutions in infinite-dimensional spaces, {\em Nonlinear Anal.} {\bf 86} (2013), 159--180.

\bibitem{MordukhovichNghia1} B. S. Mordukhovich and T. T. A. Nghia, Local monotonicity and full stability of parametric variational systems, {\em SIAM J. Optim.} {\bf 26} (2016), 1032--1059.

\bibitem{mor-roc}  B. S. Mordukhovich and R. T. Rockafellar, Second-order subdifferential calculus with applications to tilt stability in optimization, {\em SIAM J. Optim.} {\bf 22} (2012), 953--986.

\bibitem{mor-sar} B. S. Mordukhovich and M. E. Sarabi, Generalized Newton algorithms for tilt-stable minimizers in nonsmooth optimization, {\em SIAM J. Optim.} {\bf 31} (2021), 1184--1214.

\bibitem{Nikodem-Pales} K. Nikodem and  Z. P\'ales, Characterizations of inner product spaces by strongly convex functions, {\em Banach J. Math. Anal.} {\bf 5} (2011), 83--87.

\bibitem{Pen02} T. Pennanen, Local convergence of the proximal point algorithm and multiplier methods without monotonicity, {\em Math. Oper. Res.} {\bf 27} (2002), 170--191.

\bibitem{Poli} R. A. Poliquin and R. T. Rockafellar, Tilt-stability of a local minimum, {\em SIAM J. Optim.} {\bf 8} (1998), 287--299.

\bibitem{Polyak} B. T. Polyak, Existence theorems and convergence of minimizing sequences in extremum problems with restrictions, {\em Soviet Math. Dokl.} {\bf 7} (1966), 72--75.

\bibitem{Rockafellar70-Pac} R. T. Rockafellar, On the maximal monotonicity of subdifferential mappings, {\em Pacific J. Math.} {\bf 33} (1970), 209--216.

\bibitem{rtr251} R. T. Rockafellar, Variational convexity and local monotonicity of subgradient mappings, {\em Vietnam J. Math.} {\bf 47} (2019), 547--561.

\bibitem{roc} R. T. Rockafellar, Augmented Lagrangians and hidden convexity in sufficient conditions for local optimality, {\em Math. Program.} {\bf 198} (2023), 159--194.

\bibitem{r22} R. T. Rockafellar, Convergence of augmented Lagrangian methods in extensions beyond nonlinear programming, {\em Math. Program.} {\bf 199} (2023), 375--420.

\bibitem{Rockafellar98} R. T. Rockafellar and R. J-B Wets, {\em Variational Analysis}, Springer, Berlin, 1998.

\bibitem{Thibault} L. Thibault, {\em Unilateral Variational Analysis in Banach spaces},  World Scientific, Singapore, 2023.

\bibitem{ding}  S. Wang, C. Ding, Y. Zhang and X. Zhao, Strong variational sufficiency for nonlinear semidefinite programming and its applications, {\em SIAM J. Optim.}, to appear; arXiv:2210.04448 (2023).

\bibitem{Xu91} H. K. Xu, Inequalities in Banach spaces and applications, {\em Nonlinear Anal.} {\bf 16} (1991), 1127--1138.

\bibitem{Zalinescu} C. Z\u{a}linescu, {\em Convex Analysis in General Vector Spaces}, World Scientific, Singapore, 2002. 
\end{thebibliography}
\end{document}